\documentclass[11pt,leqno]{amsart}
\usepackage{amsthm,amsfonts,bm,amssymb,amsmath,color}%
\usepackage{graphicx}
\numberwithin{equation}{section}
\usepackage{mathtools,hyperref}

\mathtoolsset{showonlyrefs}

\renewcommand\d{\partial}
\renewcommand\a{\alpha}

\newcommand\s{\sigma}
\renewcommand\t{\tau}
\newcommand\R{\mathbb R}\newcommand\N{\mathbb N}\newcommand\Z{\mathbb Z}
\newcommand\C{\mathbb C}

\def\t{\tau}

\def\l{\lambda}

\def\im{\Im m\,}
\def\re{\Re e \,}
\def\epsilon{\varepsilon}
\def\e{\varepsilon}

\def\z{\zeta}

\def\t{\tau}

\def\dsp{\displaystyle}

\newcommand\br{\begin{rem}}
\newcommand\er{\end{rem}}
\newcommand\bp{\begin{pmatrix}}
\newcommand\ep{\end{pmatrix}}
\newcommand\be{\begin{equation}}
\newcommand\ee{\end{equation}}
\newcommand\ba{\begin{equation}\begin{aligned}}
\newcommand\ea{\end{aligned}\end{equation}}

\newcommand{\T}{{\mathbb T}}

\newtheorem{defi}{Definition}[section]
\newtheorem{theo}{Theorem}

\newtheorem{rem}[defi]{Remark}

\def\op{{\rm op} }

\numberwithin{equation}{section}

\begin{document}

\title{Time-delayed instabilities in complex Burgers equations}

\author{Marta Strani}
\address{Dipartimento di Matematica e Applicazioni, Universit\`a di Milano Bicocca}
\email{marta.strani@unimib.it} 

\author{Benjamin Texier}
\address{Institut de Math\'ematiques de Jussieu-Paris Rive Gauche UMR CNRS 7586, Universit\'e Paris-Diderot}
\email{benjamin.texier@imj-prg.fr}
\begin{abstract} For Burgers equations with real data and complex forcing terms, Lerner, Morimoto and Xu [{\it Instability of the Cauchy-Kovalevskaya solution for a class of non-linear systems}, Amer.~J.~Math.~2010] proved that only analytical data generate local $C^2$ solutions. The corresponding instabilities are however not observed numerically; rather, numerical simulations show an exponential growth only after a delay in time. We argue that numerical diffusion is responsible for this time delay, as we prove that for Burgers equations in the torus with small viscosity and a complex forcing, oscillating data generate solutions which grow linearly in time before growing exponentially. Numerical simulations illustrate the results.
\end{abstract}

\maketitle

\section{Introduction}

For one-dimensional Burgers equations with small viscosity, and a complex forcing term:
$$ {\rm (B)}_\e \qquad  \d_t u + u \d_x u - \e \d_x^2 u = i,$$
where $u(\e,t,x) \in \C$ depends on time $t \in \R_+,$ space $x \in \T = \R/\Z,$ and the small parameter $\e > 0,$ we prove that oscillating data generate solutions that grow linearly in time before growing exponentially. In particular, these solutions 
 behave qualitatively very much like solutions to the viscous degenerate Cauchy-Riemann equation
\be \label{cc}
 \d_t v + i t \d_x v - \e \d_x^2 v = 0,
\ee
for which the instability (of the Cauchy problem, for real data) is manifest only after a delay in time.

 This sheds new light on the strong instability result of Lerner, Morimoto, and Xu \cite{LMX} for equation ${\rm (B)}_0,$ corresponding to the inviscid case $\e = 0.$ 

\subsection{Background} \label{sec:background}   
In the inviscid case $\e = 0,$ Lerner, Morimoto and Xu proved \cite{LMX} that if a real datum generates a local $C^2$ solution of the inviscid equation ${\rm (B)}_0,$ then the datum must be analytic. This reveals a strong instability of the Cauchy problem for ${\rm (B)}_0.$

It is important to note the degeneracy in ${\rm (B)}_0:$ given a real datum $a,$ the linearized first-order operator $\d_t + a \d_x$ is hyperbolic at $t = 0;$ for arbitrarily small $t > 0,$ however, due to the complex forcing, if a solution $u(t)$ exists then presumably $\Im m \, u(t) \neq 0,$ so that the linearized operator is not hyperbolic. The degeneracy is apparent with the change of unknown $u = it + v,$ by which ${\rm (B)}_0$ becomes equation $\d_t v + (i t + v ) \d_x v = 0.$ Linearizing at $v = 0,$ we obtain the degenerate Cauchy-Riemann equation
\be \label{Btriv}
\d_t w_1 + i t \d_x w_1 = 0.
\ee
Of course, the Cauchy problem for \eqref{Btriv} is ill-posed in Sobolev spaces. The instability, however, is slower to develop than for the (non-degenerate) Cauchy-Riemann equation
\be \label{C-R}
\d_t w_2 + i \d_x w_2 = 0.
\ee
Indeed, a datum $a(x)$ will generate solutions $w_1$ to \eqref{Btriv} and $w_2$ to \eqref{C-R} which satisfy 
\be \label{sol:degCR}
 | \hat w_1(t,\xi) | = \big| \hat a(\xi) \big| e^{t^2 \xi/2}, \qquad | \hat w_2(t,\xi) | = |\hat a(\xi)| e^{t \xi},
\ee 
so that, for a given frequency $\xi$ and small $t,$ there holds $|\hat w_1(t,\xi)| \ll |\hat w_2(t,\xi)|.$  

This degeneracy is inherent to the phenomenon of transition from hyperbolicity (at $t = 0$ for ${\rm (B)}_0$) to ellipticity (at $t > 0$ for ${\rm (B)}_0$). 

The article \cite{LNT} expounds on \cite{LMX} by systematically describing the transition to ellipticity (defined as loss of hyperbolicity) for first-order systems. It is shown in \cite{LNT} that a loss of hyperbolicity implies a strong Hadamard instability, that is an instantaneous deviation estimate: some nearby data (as measured in a strong norm) may generate solutions which are instantly driven apart (as measured in a weak norm). We write ``may generate solutions", conditional, since in this setting the existence of solutions typically cannot be proved\footnote{Unless analyticity, of the data (in $x$) and the differential operator (in $u$), is assumed. For long-time Cauchy-Kovalevskaya results which prove existence together with Hadamard ill-posedness, see M\'etivier \cite{Me} in the non-degenerate case of an elliptic Cauchy problem, and Morisse \cite{Mo} in the degenerate case of a transition to ellipticity.}. Of course if some Sobolev data do not generate solutions at all, then the Cauchy problem is even more strongly ill-posed (absence of a solution operator, compared to absence of H\"older estimates for a solution operator). 

 Here we are interested in the effect of adding a small viscous term to an equation and a solution that experience a transition  from hyperbolicity to ellipticity, as we focus on ${\rm (B)}_\e,$ particularly in relation to \eqref{cc}. 
\subsection{Observations} \label{sec:numobs}

Our starting observation is that instantaneous instabilities for ${\rm (B)}_0$ are not observed numerically. Precisely, our numerical tests show for small times a linear growth of the imaginary part of the numerical solution instead of a catastropic amplification.

Our second observation is that this behavior, linear growth in time of the imaginary part for small times, is precisely the one recorded for the viscous degenerate Cauchy-Riemann equation \eqref{cc}.
Indeed, for the solution $v$ to \eqref{cc} issued from $v(\e,0,x) = a(\e,x),$ there holds, by Fourier transform and direct time integration, 
\be \label{hat:v} |\hat v(\e,t,k)| = |\hat a(\e,k)| \exp\Big( 2 \pi t k \big( \frac{t}{2} - 2 \pi \e k \big)\Big), \qquad k \in \Z,\ee
so that $\hat v(\cdot,k)$ grows exponentially only for $t \geq 4 \pi \e k,$ and only if $k$ belongs to the support of $\hat a.$ 

For instance, if the initial datum is highly oscillating: $a(x) = \sin(2 \pi x/\e),$ then the amplification occurs only after $t = 4 \pi.$ And if the smallest positive frequency in the datum is $k_0 \in \Z_+^*,$ as in $a(x/\e) = \sin(2 \pi k_0 x/\e),$ then the amplification for \eqref{cc} occurs only after $t = 4 \pi k_0.$

We define the {\it transition time} as the smallest positive time for which the {\it symbol has negative real part.} We compute the symbol associated with equation \eqref{cc} by changing $\d_x$ into $2 \pi i \xi$ (Fourier transform), so that the symbol here is   
%
$$ i t (2 \pi i\xi) - \e (2 \pi i \xi)^2 = - (2\pi \xi) \big( t - 2 \pi \e \xi\big).$$
Real negative values occur for $t \geq 2 \pi \e \xi.$ For instance, for the datum $a(\e,x) = \sin(2 \pi x k_0),$ with $k_0 \in \Z^*,$ we consider frequencies $\xi \in k_0 \Z,$ and the smallest $t$ for which the symbol has real negative eigenvalues is $2 \pi \e k_0.$ This is the {\it transition time} (for this equation and this datum). For the datum $a(\e,x) = \sin(2 \pi k_0x/\e),$ the transition time is $2 \pi k_0.$ 

By {\it amplification time} we mean smallest possible time for which an adequate measure of the solution is greater than the same measure of the datum. For instance, in the case of datum $a(\e,x) = \sin(2 \pi k_0 x/\e),$ the representation \eqref{hat:v} for the solution $v$ to \eqref{hat:v} is non-trivial only if $k \in k_0\Z/\e.$ 
In particular, the smallest possible $t$ for which $|\hat v(\e,t,k)| \geq |\hat a(\e,k)|$ {\it for some $k$} is $t = 4 \pi k_0,$ corresponding to $k = k_0/\e.$ The amplification time is $4 \pi k_0,$ greater than the transition time.

In particular, the amplification occurs {\it only after} the transition. This is particularly meaningful in view of our results stated in Section \ref{sec:results} below.  

\medskip

Based on these observations, our guess was that numerical diffusion was responsible for the defect in instantaneous amplification in the simulations. By numerical diffusion, we mean the fact that the standard Lax-Friedrichs scheme that we used for the simulation of ${\rm (B)}_0,$ with time step $\s$ and length interval $h,$ is consistent at order $1$ with ${\rm (B)}_0,$ but consistent at order 2 with ${\rm (B)}_\e$ with $\e = h^2/(2 \s).$ 

 The analogy with \eqref{cc} then suggested that exponential amplification for the numerical solutions to ${\rm (B)}_0$ would occur for ulterior times. This was confirmed by numerical simulations on longer time intervals. 

\begin{figure}
\includegraphics[scale=.9]{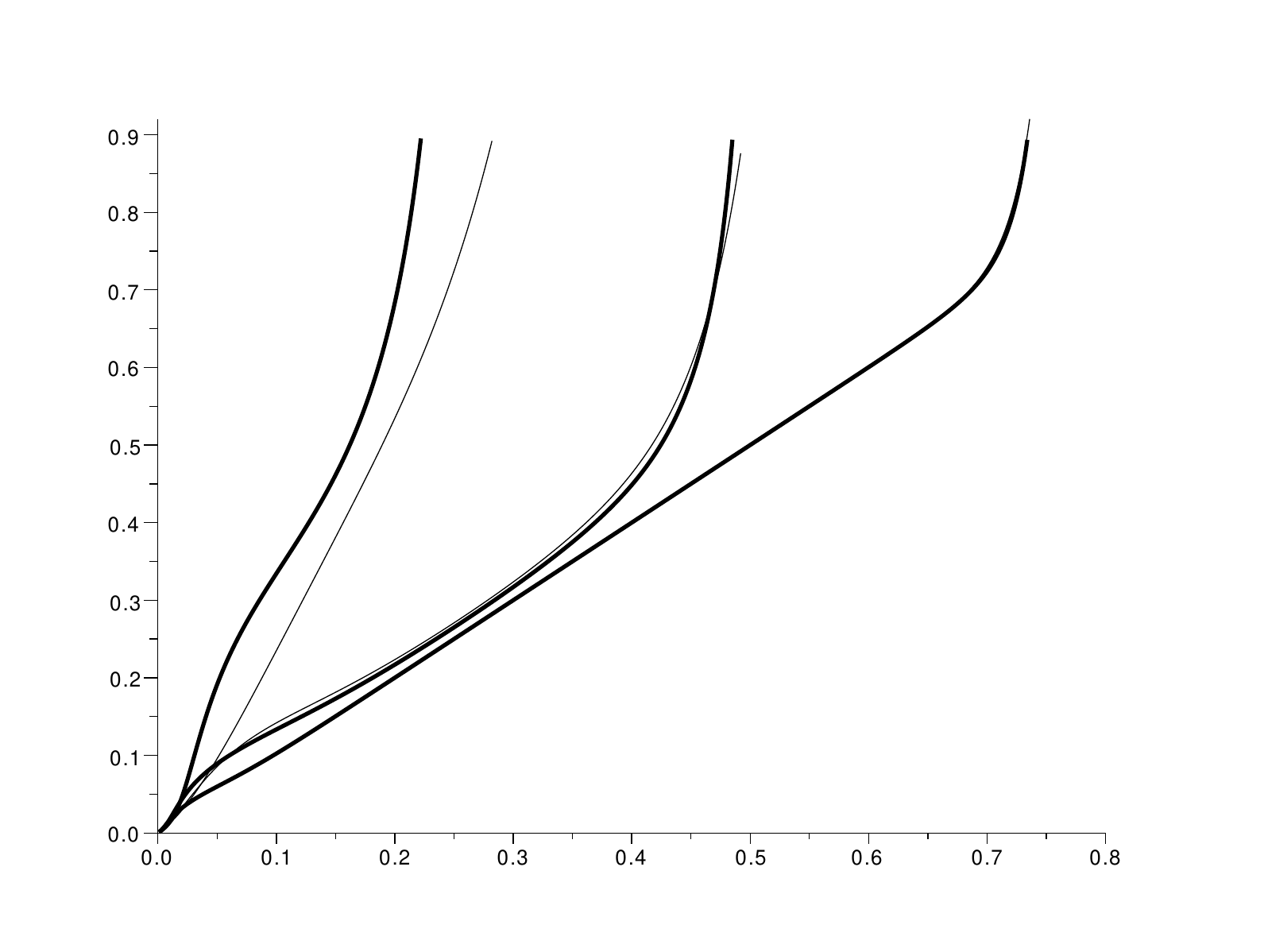}
\caption{Numerical and ``linearized" solutions for $N = 8, 16$ and $24.$}
\label{fig-1}
\end{figure}

 Numerical solutions to ${\rm (B)}_\e$ are pictured on Figure \ref{fig-1}. We show tests for the initial data $a(x) = \sin(N \cdot 2 \pi x),$ with $N \in \{ 8, 16, 24 \}.$ We recorded the maximum value over the spatial grid of the imaginary parts of the numerical solutions every 20 time steps. The thick lines are the graphs of these maxima as functions of time. The solution that seems to blow up around $t = .2$ corresponds to $N= 8.$ The one that seems to blow up around $t=.45$ corresponds to $N = 16.$ The rightmost one, seemingly blowing up around $t = .7,$ corresponds to $N = 24.$ By comparison, we picture on the same graph the corresponding ``linearized" solutions, in thin lines. These solutions are translates of solutions to \eqref{cc}. Indeed, if $u$ solves ${\rm (B)}_\e,$ then $v := u - i t$ solves $\d_t + (i t + v) \d_x v - \e \d_x^2 v = 0,$ and linearizing at $v = 0$ we find \eqref{cc}. The solution to \eqref{cc} issued from $a(x) = \sin(N \cdot 2 \pi x)$ is 
  $$ v(\e,t,x) = \frac{1}{2i} e^{t^2 N/2 - \e t N^2} e^{2 i \pi N x} - \frac{1}{2i} e^{-t^2 N/2 - \e t N^2} e^{-2 i \pi N x},$$
so that
$$ \max_{x \in \T} \im v(\e,t,x) = \frac{1}{2} \Big( e^{t^2 N/2 - \e t N^2}  - e^{-t^2 N/2 - \e t N^2} \Big).$$
The thin lines on Figure \ref{fig-1} are the graphs of $t \to t + \max_{\T} \im v(\e,t,\cdot),$ with $v$ as above. For $N = 16$ we find good agreement between the numerical and the linearized solution. For $N = 24,$ on the scales of the picture the graphs are essentially indistinguishable. This indicates that the approximation of ${\rm (B)}_\e$ by \eqref{cc} seems quite accurate for oscillating data, even if the amplitude of the data is $O(1).$ 

Here the ratio $\e = h^2/(2 \s)$ is set equal to $2.5 \cdot 10^{-3},$ with the spatial step $h = 5 \cdot 10^{-4}$ and the temporal step $\s = 5 \cdot 10^{-5}.$

On Figure \ref{fig-1}, the rightmost graph, corresponding to $N = 24,$ shows a linear behavior in time, followed by a quick amplification.  The computing time for this simulation is about 250 minutes on a laptop computer. The amplification is apparent only about $t = .7,$ well after 150 minutes of computing time.

Next we examined how the amplification time varied according to the initial oscillations and the viscosity. According to \eqref{hat:v}, the amplification time for equation \eqref{cc} is $t_{\rm cc} = 4 \pi \e N,$ where $N$ is the smallest non-zero mode in the real initial datum. In particular, it has a linear dependence both in $\e$ and in the smallest mode in the datum. We compared this linear amplification time with $t_{\rm f},$ the final computing time in our simulations, defined as follows.

In the main time-loop of our numerical scheme, we imposed the CFL condition
 \be \label{cfl} \frac{\s}{h} = 0.1 < \frac{0.4}{\max(\max_j |\re u^n(j)|, \max |\im u^n(j)|)},\ee
 where $u^n$ represents the numerical solution to ${\rm (B)}_\e$ at step $n.$ In particular, the simulations stopped when condition \eqref{cfl} broke down, corresponding to the smallest positive time for which the $L^\infty$ norm of the discrete solution is greater or equal to 4. The final computing time $t_{\rm f}$ is the largest computing time before \eqref{cfl} breaks down.

\begin{figure}
\includegraphics[scale=.6]{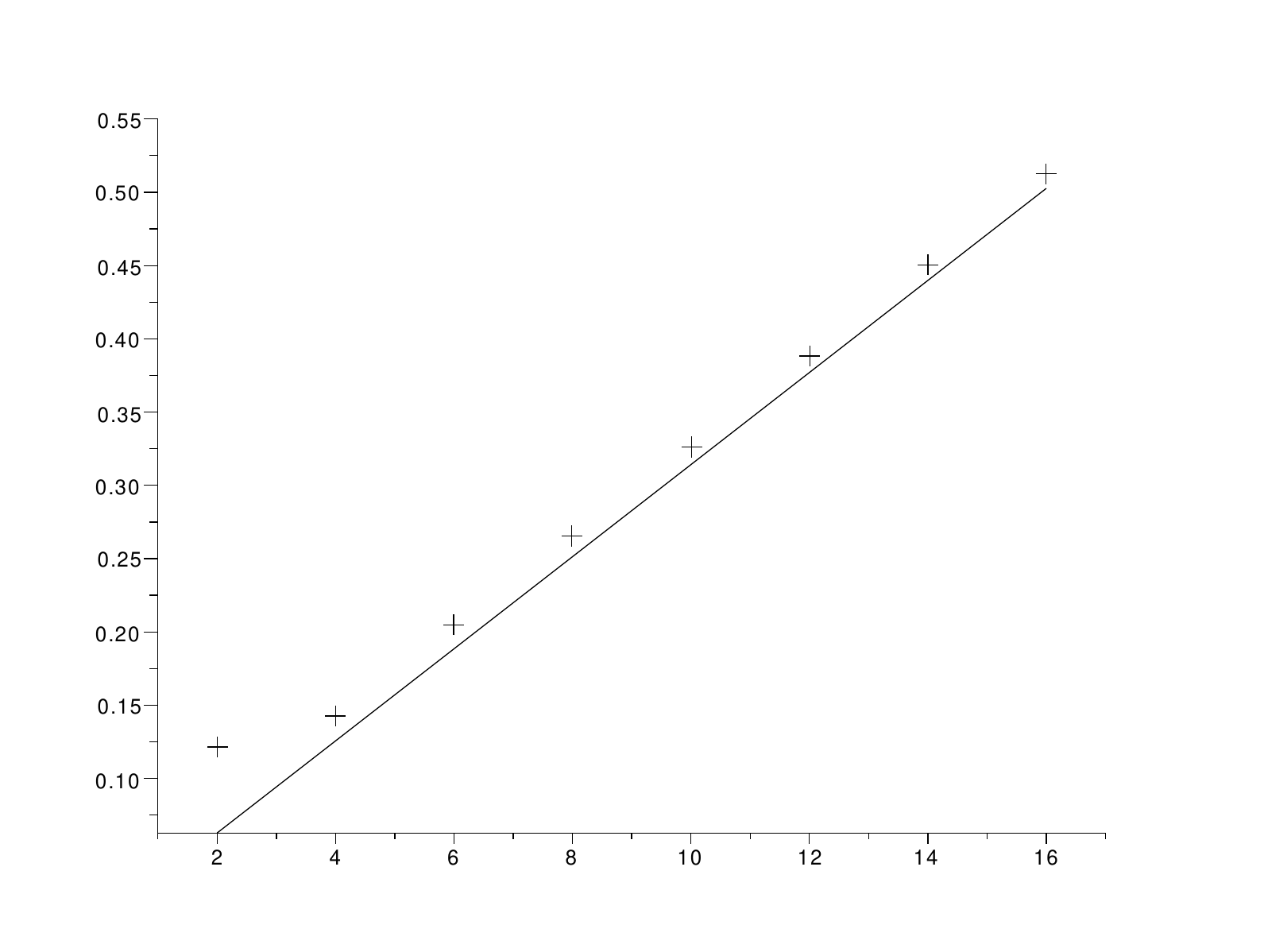}
\caption{Amplification times for different values of the initial frequency.}
\label{fig-4}
\end{figure}

On Figure \ref{fig-4} we compare $t_{\rm cc} = 4 \pi \e N$ (drawn as full line) with $t_{\rm f}$ (crosses), for the initial data $a(x) = \sin(N \cdot 2 \pi x),$ for $N$ taking even values from $2$ to $16.$ The values of $h$ and $\s,$ hence $\e,$ are as in Figure \ref{fig-1}. We see that for $N$ not too small, there is good agreement between $t_{\rm cc}$ and $t_{\rm f}.$ It is not surprising to find $t_{\rm f}$ to be greater than $t_{\rm cc},$ since when the maximum of the numerical solution equals $4,$ the amplification has already started for some time. On Figure \ref{fig-1}, the graphs are truncated before the final computing time for better legibility.

\begin{figure}
\includegraphics[scale=.6]{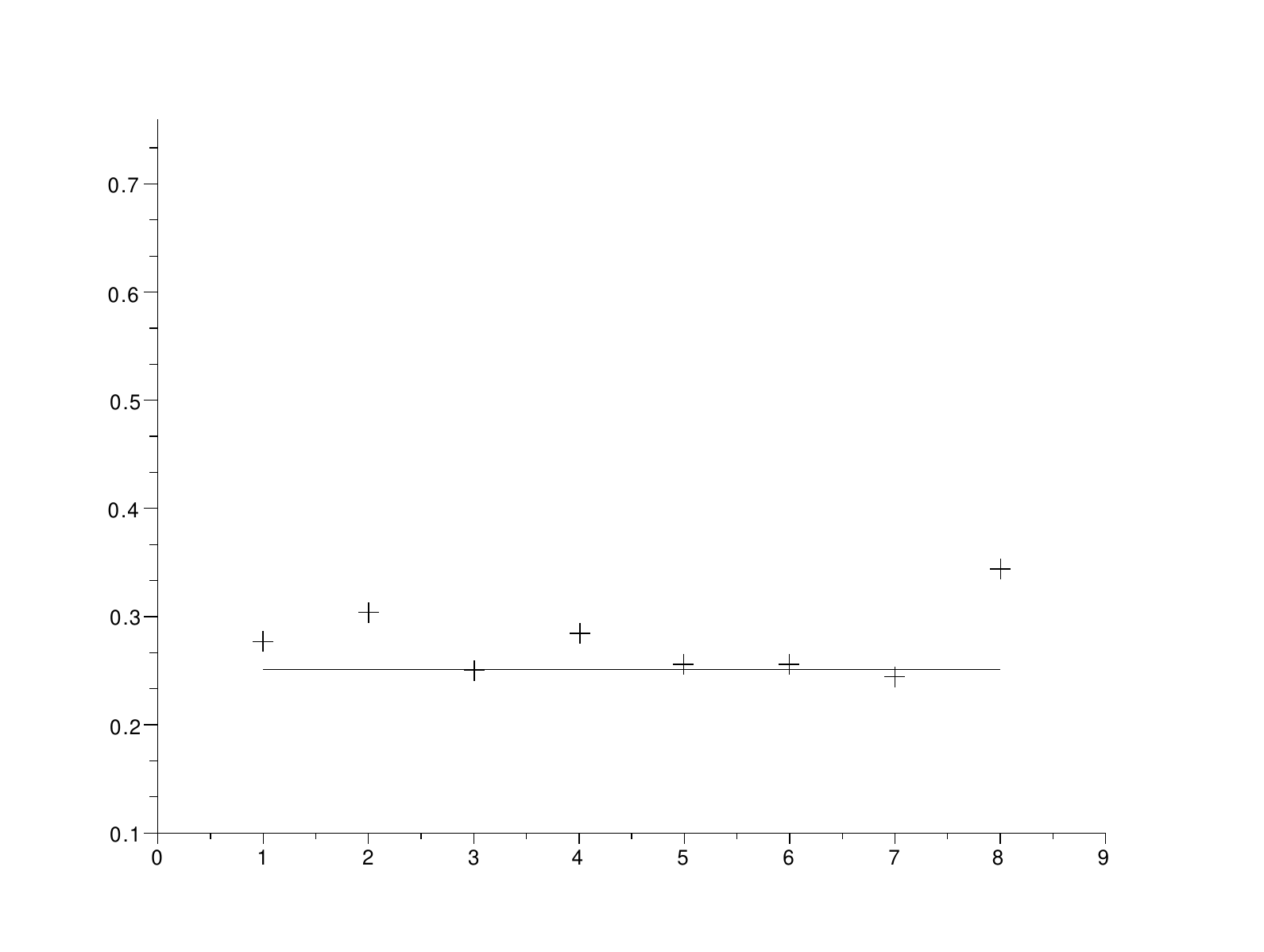}
\caption{Amplification times for different initial data.}
\label{fig-5}
\end{figure}

On Figure \ref{fig-5}, we show tests for the dependence of the amplification time on the smallest non-zero mode in the datum. We compare $t_{\rm cc}$ (full line) with $t_{\rm f}$ (crosses), for the same values of $h$ and $\s$ as above, and for 8 different initial data, of the form $\sum_{1 \leq j \leq 3} a_j \sin(N_j \cdot 2 \pi x),$ with specific values as follows, the $x$-coordinate on Figure \ref{fig-5} corresponding to the case number:
\begin{itemize} 
\item case 1 : $(N_1,N_2,N_3) = (4,0,0),$ $(a_1,a_2,a_3) = (1,0,0);$ this is a test case with only one initial oscillation;
\item case i, with $2 \leq i \leq 5:$ $(N_1,N_2,N_3) = (4, 4 + 2(i-1),0),$ $(a_1,a_2,a_3) \equiv (1,1,0);$
\item cases 6,7,8: $(N_1,N_2,N_3) \equiv (4,6,8),$ with coefficients equal to $(1,1,1)$ then $(1,2,1)$ and finally $(1,2,-1).$
\end{itemize}
The straight line is the linear amplification time $t_{\rm cc} = 4 \cdot 4 \pi \e,$ common to all cases. We observe that the numerical values of $t_{\rm f}$ are relatively close to the linear value. Of course the scale on the $y$-axis here plays a crucial role in conveying the impression that the approximation by \eqref{cc} makes sense. The 8 cases tested on Figures \ref{fig-5} involve oscillations $\sin(2 \pi N x),$ with $N$ ranging from $4$ to $12.$ Therefore we chose extremal values on the $y$-axis corresponding to the extremal amplification times associated with these frequencies, that is values ranging from $t_{\rm cc}(N = 4) = .251$ to $t_{\rm cc}(N = 12) = .754.$

\begin{figure}
\includegraphics[scale=.6]{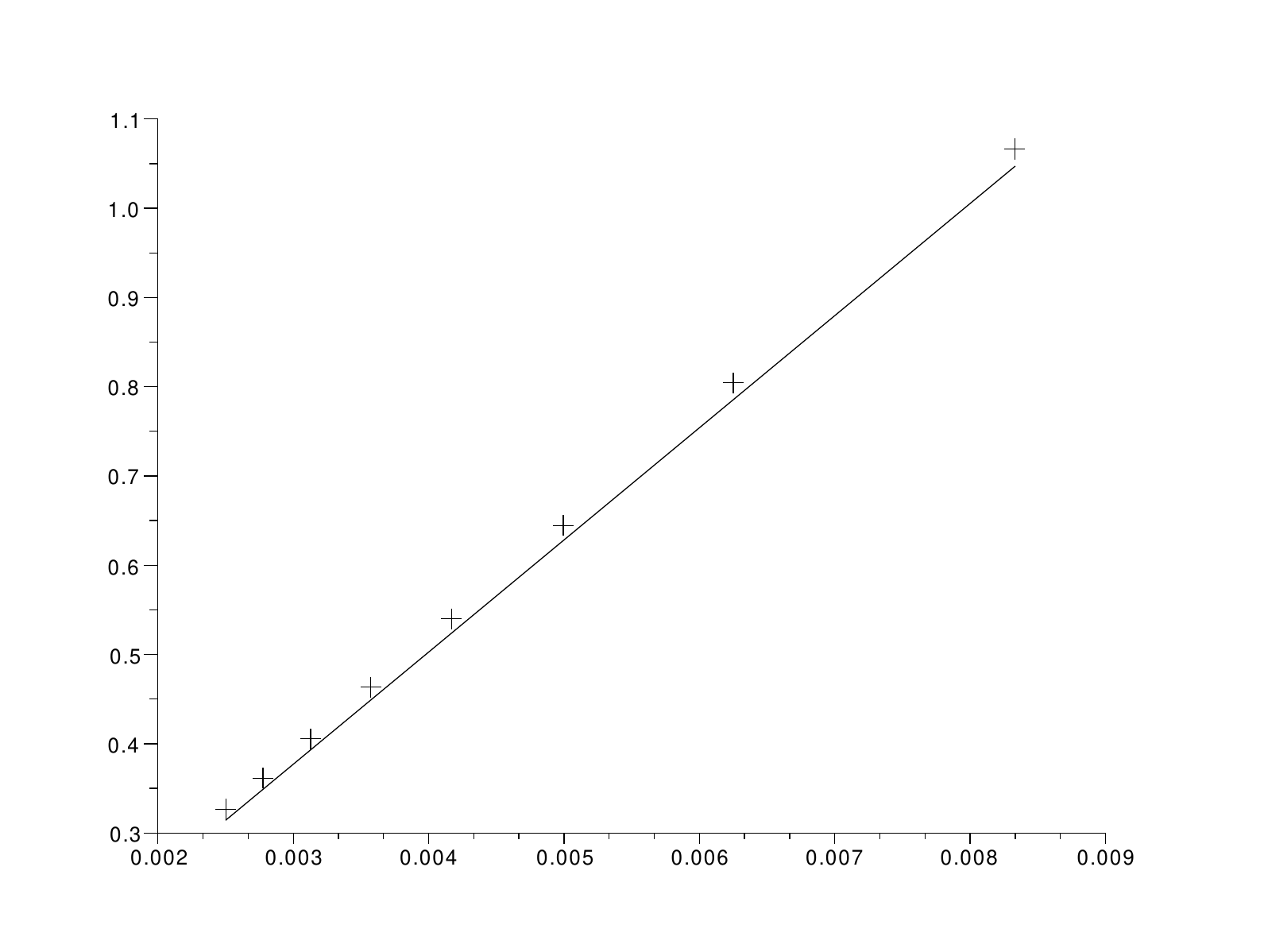}
\caption{Amplification times depending on viscosity.}
\label{fig-6}
\end{figure}

We finally tested the dependence of the final computing time on the viscosity. The results are pictured on Figure \ref{fig-6}, with $t_{\rm cc} = 4 \pi \e N$ drawn as a full line, and $t_{\rm f}$ as crosses. Here the datum is $a(x) = \sin(10\cdot 2 \pi x).$ We tested values of the number $J$ of points on the spatial grid, from $J = 2000$ to $J = 600,$ with increments of 200. The associated spatial step is $h = 1/J.$ We fixed $\s = h/10,$ so that $\e = h^2/(2 \s) = 5/J.$ The corresponding values of $\e$ range from $2.5 \cdot 10^{-3}$ (for $J = 2000$) to $8.33 \cdot 10^{-3}$ (for $J = 600$). There is an excellent agreement with the linear approximation.  
 
\subsection{Results} \label{sec:results}

In a first step (Section \ref{sec:crude}), we show by elementary a priori bounds that non-oscillating data $a(x)$ generate solutions that are bounded in short time $t < 2 \pi \e.$ Then we consider highly-oscillating data, with characteristic frequencies $O(1/\e),$ and show that these generate solutions that are linearly growing up to the transition time $2 \pi k_0,$ where $k_0$ is the leading mode in the datum:
\begin{theo} \label{th1} Given $k_0 \in \N^*,$ given $0 < T < 2 \pi k_0$ and an integer $s \geq 2,$ given $a \in H^{s+1}(\T),$ if $\e$ and $|a - \bar a|_{H^{s+1}}$ are small enough the datum 
$$ u(\e,0,x) = a\left(\frac{k_0 x}{\e}\right)$$
generates a unique solution $u \in  C^0([0, T], H^{s}(\T))$ to ${\rm (B)}_\e,$ and there holds
$$ u(\e,t,x) = i t + v\left(\e, \, \frac{t}{\e}, \,  \frac{x}{\e}\,\right),$$
with bounds  
\be \label{growth:est} \begin{aligned} |\Re e \, v(t) - \bar a |_{H^s} & \lesssim \big (1 + \e t C(T) \big) |a(k_0 \cdot) - \bar a |_{H^{s+1}}, \\ |\Im m \, v(t) |_{H^s} & \lesssim \e t \big(1 + \e t C(T) \big) |a(k_0\cdot) - \bar a|_{H^{s+1}}, \end{aligned}
\ee
where $0 < C(T) \to \infty$ as $T \to 2 \pi k_0.$ \end{theo}


In the statement of Theorem \ref{th1}, $\bar a$ denotes the mean mode of the initial datum: $\dsp{\bar a := \frac{1}{2\pi} \int_{\T} a(x) \, dx}.$ The implicit multiplicative constants in \eqref{growth:est} depend only on $s$ and $k_0.$ Theorem \ref{th1} is proved in Section \ref{sec:proof1}. 

\medskip

Based on the numerical tests, we expect the solutions of Theorem \ref{th1} to grow exponentially for $t > 2 \pi k_0.$ This is something we can prove for small data: 
\begin{theo} \label{th2} Given $k_0 \in \N^*,$ given $T > 0,$ $\a > 1/3,$ and $a \in H^1(\T)$ such that $\dsp{\int_{\T} e^{- 2 i \pi x} a(x) \, dx \neq 0},$ if $\e$ is small enough the datum 
  $$ u(\e,0,x)  = \e^{\a} a\left(\frac{k_0 x}{\e}\right)$$
generates a unique solution $u \in C^0([0, 4 \pi k_0 + \e T], H^1(\T))$ to ${\rm (B)}_\e,$ and there holds
 $$ \left| \int_{\T} e^{-2 i \pi k_0 x/\e} u(\e,t,x) \, dx \right| \geq \frac{1}{2} e^{\pi k_0 ( t - 2 \pi k_0)^2/\e - 4 (\pi k_0)^3/\e} \left| \int_{\T} e^{- 2 i \pi k_0 x/\e} u(\e,0,x) \, dx \right|.$$
   \end{theo} 

Theorem \ref{th2} is proved in Section \ref{sec:growth}. The general case of $O(1)$-amplitude data seems out of reach; we explain why in Section \ref{sec:open}. Note that the time-exponential lower bound in Theorem \ref{th2} starts to grow at $t = 2 \pi k_0,$ the {\it transition time} described in the introduction, and reaches the value $1$ at $t = 4 \pi k_0,$ the {\it amplification time} described in the introduction. Our description of the growth stops shortly after the amplification time.

\subsection{Comments} \label{sec:comments} {\it On Theorem {\rm \ref{th1}} and its proof.} The proof is based on the construction of an approximate solution to ${\rm (B)}_\e.$ This procedure typically has a cost in terms of regularity. This shows in Theorem \ref{th1}, as we need $a - \bar a$ to belong to $H^{s+1},$ while we perform estimates only in $H^s.$ The proof uses only elementary inequalities, such as Poincar\'e-Wirtinger and Kato-Ponce. We note in Section \ref{sec:crude} that existence and bounds up to the transition time are easy to derive. Our point in Theorem \ref{th1} is to obtain precise growth estimates for the imaginary part of the solution; in this view compare the crude bound \eqref{crude:bound} obtained at the end of our basic analysis in Section \ref{sec:crude} with estimate \eqref{growth:est}. 

 Theorem \ref{th1} is illustrated on Figure \ref{fig-1} by the behavior of the numerical solutions for $N = 16$ and $N = 24$ before their respective amplification times, and also by Figure \ref{fig-5}: the amplification time depends chiefly on the smallest non-zero mode in the datum. Here the smallest mode is $k_0/\e,$ corresponding to a linear amplification time $4 \pi k_0,$ as discussed in Section \ref{sec:numobs}.
 
This slow-then-fast dynamics is reminiscent of the phenomenon of {\it metastability}, in which a speed of convergence depends singularly on the viscosity, typically as $e^{-1/\e}$\cite{KK,meta}. These two phenomena, metastability and
time-delay in the instabilities, have in common the fact that the solutions
exhibit a certain stable (observable)
behavior for an $O(1)$ time interval before they converge to the asymptotic
limit in the case of a metastable behavior, or experience an exponential
growth in time in the case of a time-delay in the instability.

\smallskip

{\it On Theorem {\rm \ref{th2}} and its proof.} For small solutions, the approximation by the linear constant-coefficient equation \eqref{cc} certainly makes sense at least for small time. Theorem \ref{th2} shows that qualitatively, the approximation is actually valid up to and a little bit beyond the amplification time $t = 4 \pi k_0,$ which is not clear at all at first sight.  

For \eqref{cc}, the mode-by-mode analysis \eqref{hat:v} shows that the Fourier modes $k/\e < 0$ decay exponentially, while Fourier modes $k/\e > 0$ decay until $t = 2 \pi k$ (transition time), recover their initial value at $t = 4 \pi k$ (amplification time), then grow beyond their initial value, at a rate $\sim e^{\pi t k^2},$ for large $t.$ 

The proof of Theorem \ref{th2} is based on this simple description. In particular, it is conducted entirely on the Fourier side. The smallness assumption allows to handle the nonlinear convective terms as small sources. We perform elementary estimates based on the exact description of the linear solution operator in Fourier, as seen on \eqref{hat:v}. The difficulty comes from error function estimates that arise in bounds for Gaussian time integrals.

\subsection{Open problem} \label{sec:open} Our numerical observations lead us to conjecture that%
$$ \begin{array}{cl} {\rm (C)} & \begin{array}{ll} \mbox{\it Data with amplitude $O(1)$ and frequencies $O(1/\e),$ of the form $a(x/\e),$} \\ \mbox{\it generate solutions to ${\rm (B)}_\e$ which are bounded in time $O(1)$ before} \\ \mbox{\it exhibiting an exponential growth in time.}\end{array} \end{array}$$

We managed to prove the part about linear growth (Theorem \ref{th1}), but we could prove exponential amplification only for small data (Theorem \ref{th2}). We describe here in some detail the difficulties that arise in a sketch of proof of conjecture (C). 

Let $u(\e,0,x) = a(x/\e)$ be given in $H^{s+1}(\T).$ By Theorem \ref{th1}, we know that the corresponding solution $u$ is defined and bounded up to time $T < 2 \pi.$ At $t = 2\pi,$ there is a change in behavior in the associated linear constant-coefficient equation \eqref{cc}. Indeed, as we saw indirectly in \eqref{hat:v}, the symbol $\l(\e,t,\xi) = \e (2 \pi \xi)^2 - 2 \pi t \xi $ of the differential operator in \eqref{cc} satisfies $\l > 0$ for $t < 2 \pi$ (damping) and $\l < 0$ for $t > 2 \pi$ (amplification).

Assume that $u$ is defined and bounded up to and a little bit beyond the transition time $t = 2 \pi.$ We posit $u = i t + \bar a + v,$ so that $v$ solves
\be \label{eq:v:intro}
 \d_t v + \frac{1}{\e} \big( i t + \bar a + v\big) (\e \d_x) v - \frac{1}{\e} (\e \d_x)^2 v = 0.
\ee
Since $v$ oscillates rapidly in $x$ (like its datum, presumably), it makes sense to write \eqref{eq:v:intro} in terms of $\e \d_x$ derivatives, rather than $\d_x$ derivatives. The symbol $\l(\e,t,x,\xi)$ of the differential operator in \eqref{eq:v:intro} is
$$ \l(\e,t,x,\xi) = \big( i t + \bar a + v(\e,t,x) \big) 2 i \pi \xi + (2 \pi \xi)^2,$$
so that \eqref{eq:v:intro} takes the form 
\be \label{eq:v:intro:bis}
 \d_t v + \frac{1}{\e} \op_\e(\l) v = 0, \qquad \op_\e(\l) v := \int_{\R^d} e^{2 i \pi x \cdot \xi} \l(\e,t,x, \e \xi) \hat v(\xi) \, d\xi.
\ee
Now if locally in $(t,x,\xi),$ there holds $\Re e \, \l < 0,$ we expect a fast exponential growth for \eqref{eq:v:intro:bis}.  The mode $\xi = 1$ should be amplified first, so we will look at $\l$ near $\xi = 1.$ There holds
$$ \re \l(\e,t,x,1) = 2 \pi \big( 2 \pi - t - \im v(\e,t,x) \big).$$
Since $\bar v \equiv 0,$ we know that $\im v$ is not everywhere negative. In particular, for any $t > 2 \pi,$ we can find $x = x(t)$ such that $\re \l(\e,t,x(t),1) < 0.$ We even expect $\re \l$ to take negative values {\it before} $t = 2\pi,$ since the maximum of $\im v$ is nonnegative. That is, if we let 
$$ t_\star := \inf \Big\{ 0 < t \leq 2 \pi, \quad 2 \pi - t - \max_{\T} \im v(\e,t) < 0 \Big\},$$
then $t_\star < 2 \pi$ (unless $\im v \equiv 0$ in a neighborhood of $t = 2 \pi,$ unlikely), and $t_\star$ appears as the transition time. At this stage, we wait a little bit longer for the transition to happen. That is, we do not handle the difficult case of a transition to ellipticity (central to \cite{LNT,Mo}), and instead observe the solution at some ulterior $t_0 \in (t_\star, 2\pi),$ when presumably $\re \l$ is locally negative and bounded away from zero. Let $x_0$ such that $\max_{\T} \im v(\e,t_0) = \im v(\e,t_0,x_0).$ We can localize around $x_0 = 0$ (by multiplication by a cut-off $\theta$ in $x$), then identify the solution to the truncated equation with a compactly supported function defined on $\R.$ Then we can localize in frequency around $\xi_0 = 1,$ by multiplication of the equation to the left by $\op_\e(\chi),$ where $\chi$ is a cut-off in $\xi.$ The cut-off procedures produce error terms, in the form of commutators, which are smaller by a factor $\e$ since the quantization in \eqref{eq:v:intro:bis} is semiclassical\footnote{By composition of pseudo-differential operators in semiclassical quantization. This is already observed for differential operators, in equality $a(x) \e \d_x \circ b(x) \e \d_x = a(x) b(x) (\e \d_x)^2 + \e a(x) b'(x) \e \d_x,$ where $\e a b' (\e \d_x)$ appears a a small remainder. For a general pseudo-differential composition result, see for instance Theorems 1.1.5 and 1.1.20 in \cite{Le}.}. Thus we arrive at 
\be \label{eq:trunc}
\d_t w + \frac{1}{\e} \op_\e(\mu) w = f,
\ee
where $f$ is a sum of commutators and is bounded in $\e,$ and $\mu$ has the form
$$ \mu(\e,t,x,\xi) = \theta^\sharp(x) \chi^\sharp(\xi) \l(\e,t,x,\xi),$$
where $\theta^\sharp$ is an extension of the original spatial cut-off $\theta,$ and similarly $\chi^\sharp$ is an extension of the original frequency cut-off $\chi.$ If the cut-offs are carefully chosen, then by choice of $(t_0,x_0),$ there holds $\Re e \, \mu \leq - \mu_0 < 0,$ and the growth of $w$ follows by an application of G\r{a}rding's inequality (see for instance Theorem 4.32 in \cite{Z}). We conclude that $v$ must grow as well, and in the end that it was not reasonable to assume uniform bounds for $u$ much beyond the transition time (presumably, uniform bounds break down just after the amplification time $t = 4 \pi$). 

The difficulty that we overlooked in this discussion is the fast dependence of $v,$ hence of the symbol $\l,$ on the spatial variable $x.$ In particular, we can certainly find $(t_0,x_0)$ such that $\re \l(\e,t_0,x_0,1) < 0,$ but this inequality is unlikely to hold in a $O(1)$-neighborhood of $x_0.$ Thus we would need to use an $\e$-dependent cut-off $\theta,$ typically of the form $\theta((x - x_0)/\e).$ Such a cut-off generates commutator that are {\it not} bounded in $\e.$ Thus in \eqref{eq:trunc} we will actually obtain a source $(1/\e) f,$ with a bounded $f.$ In this setting, we would need more information on $v,$ and $f,$ in order to conclude by an application of G\r{a}rding's inequality. 

 
 A similar open problem (amplification proof in the presence of wildly varying coefficients) was raised in Section 1.8 of \cite{em4}.
 
 We note another difficulty, not unrelated to the small spatial scale of variation of $v.$ It is that we need to have {\it some} information on $v$ at $t_0$ in order to prove growth. Indeed, a space-frequency truncated $v$ at $t = t_0$ will serve as a new initial datum for \eqref{eq:trunc} (with $(1/\e) f$ instead of $f$). We will need in particular $\op_\e(\chi) (\theta v)$ not to be too small at $t = t_0.$ 
 
 In conclusion, the difficulty here is that we are looking at the behavior of fast-oscillating solutions to nonlinear equations in time $O(1),$ well after nonlinear effects have started taking place. 

\section{Basic analysis: non-oscillating data and short-time bounds} \label{sec:crude} 

We prove here that a non-oscillating datum $u(0,x) = a(x) \in H^s(\T)$ generates a unique solution $u$ up to time $T,$ for any $T < 2 \pi \e,$ if $|a - \bar a|_{H^s}$ is small enough and $s \geq 1.$ For further use we note that there holds for all $u \in H^1(\T)$ the Poincar\'e-Wirtinger inequality
 \be \label{pw}
  |u - \bar u|_{L^2} \leq \frac{1}{2\pi} |\d_x u|_{L^2}, \quad \bar u := \frac{1}{2\pi} \int_\T u \, dx.
 \ee
 We look for the solution $u$ to ${\rm (B)}_\e$ in the form
$$ %
 u(t,x) = i t + \bar a + v(t,x).%
$$ %
Then $v$ solves 
$$ %
\d_t v + (i t + \bar a + v) \d_x v -  \e \d_x^2 v  = 0, \qquad \bar v \equiv 0.
$$ %
Differentiating, we obtain with $v_\a := \d_x^\a v,$ for $0 \leq \a \leq s,$ with $1 \leq s:$ 
$$ %
\frac{1}{2} \d_t \big( |v_\a|_{L^2}^2 \big) + \e |\d_x v_\a|_{L^2}^2 + \Re e \, \Big( \, (i t + \bar a) (\d_x v_\a, v_\a)_{L^2} \, \Big) + \Re e \, (\d_x^\a (v \d_x v), v_\a)_{L^2} = 0,
$$ %
where the $L^2$ scalar product is hermitian: $\dsp{(f,g)_{L^2} := \int_{\T} \bar f \, g .}$ In particular, since $\bar a \in \R,$ there holds $\Re e \, (\bar a \d_x v_\a, v_\a)_{L^2} = 0,$ so that
$$ \Big| \Re e \, \Big( \, (i t + \bar a) (\d_x v_\a, v_\a)_{L^2} \, \Big) \Big| = \Big| - t \Im m \, (\d_x v_\a, v_\a) \Big| \leq t |\d_x v_\a|_{L^2} |v_\a|_{L^2}.$$
Besides,
$$ %
 \Big| \Re e \, (\d_x^\a (v \d_x v), v_\a)_{L^2} \Big| \leq \frac{1}{2} |v^2|_{H^{\a+1}} |v_\a|_{L^2}.
$$ %
By the Moser product inequality\footnote{In the form $|w_1 w_2|_{H^s} \leq C_s (|w_1|_{L^\infty} |w_2|_{H^s} + |w_2|_{L^\infty} |w_1|_{H^s}),$ for $s \geq 0,$ some $C_s > 0,$ all $w_i \in L^\infty \cap H^s.$},
\be \label{moser} |v^2|_{H^{\a+1}} \leq C_s |v|_{L^\infty} |v|_{H^{\a+1}}.\ee
We will henceforth use $C_s$ as a generic notation for positive constants depending only on $s.$ Then by Poincar\'e-Wirtinger \eqref{pw}, since $v$ has zero mean value, $|v|_{H^{\a+1}} \leq (1 + 1/(2\pi)) |\d_x v|_{H^s}.$ 
Gathering the above bounds, we obtain
$$ %
 \frac{1}{2} \d_t \big( |v_\a|_{L^2}^2 \big) + \e |\d_x v_\a|_{L^2}^2 \leq t |\d_x v_\a|_{L^2} |v_\a|_{L^2} + C_s |v|_{L^\infty} |\d_x v|_{H^s} |v_\a|_{L^2}.
$$ %
Using Poincar\'e-Wirtinger again and then summing\footnote{Here multiplicative constants do matter, since our goal is to reach precisely the transition time $2 \pi \e.$ In particular, we define the $H^s$ norm as a euclidian norm: $|w|_{H^s}^2 = \sum_{\a \leq s} |\d_x^\a w|_{L^2}^2.$} over $\a \leq s$, this gives
\be \label{bd:0}
 \frac{1}{2} \d_t \big(|v(t)|_{H^s}^2\big) + \e |\d_x v|_{H^s}^2 \leq \left( \frac{t}{2\pi} + C_s |v|_{L^\infty} \right) |\d_x v|_{H^s}^2.
 \ee
In particular, so long as
\be \label{cond:0}
 \frac{t}{2\pi} + C_s |v|_{L^\infty} \leq \e,
\ee
there holds 
\be \label{up:1}
 |v(t)|_{H^s} \leq |a - \bar a|_{H^s}.
  \ee
Thus for any $0 < T < 2 \pi \e,$ the uniform bound \eqref{up:1} holds over all of $[0, T],$ provided that $| a - \bar a|_{H^s} \leq C_s(2 \pi \e - T).$ Standard arguments can be used to convert this uniform a priori estimate into an existence and uniqueness result in $L^\infty([0, T], H^s(\T)),$ for any integer $s \geq 1.$

\section{Proof of Theorem \ref{th1}: oscillating data and uniform bounds in time $O(1)$} \label{sec:proof1}

We consider here oscillating data, of the form $\dsp{u(\e,0,x) = a\left(\frac{k_0 x}{\e}\right)},$ where $k_0 \in \N.$ We look for the solution $u$ to ${\rm (B)}_\e$ issued from this datum in the form 
\be \label{scales}
 u(\e,t,x) = i t + \bar a + v\left(\frac{t}{\e}, \frac{x}{\e}\right),
\ee 
corresponding to a hyperbolic change of scales. Then $v$ solves 
\be \label{eq:vtilde} %
\d_t v + (i \e t + \bar a + v) \d_x v -  \d_x^2 v  = 0, \qquad \bar v \equiv 0, \qquad v(\e,0,x) = a(k_0 x) - \bar a.
\ee %
We observe that the datum can be expanded in Fourier series
$$ v(\e,0,x) = \sum_{ k \in \Z^* } e^{2 i \pi k \cdot x k_0} a_k,$$
where $(a_k)$ are the Fourier coefficients of $a.$ The convective nonlinearity in \eqref{eq:vtilde} produces only harmonics of $k_0.$ In particular, over its interval of existence, the solution $v$ to \eqref{eq:vtilde} enjoys a Fourier expansion
$$ v(\e,t,x) = \sum_{k \in \Z^*} e^{2 i \pi k \cdot  x k_0} v_k(\e,t),$$
and the Poincar\'e-Wirtinger inequality \eqref{pw} takes the particular form
\be \label{pw2}
 |v|_{L^2} \leq \frac{1}{2 \pi k_0} |\d_x v|_{L^2}.
\ee 
At this point we can reproduce the arguments of Section \ref{sec:crude} and prove existence up to time $T < 2 \pi k_0.$ Indeed, by the same arguments as above, taking into account the scaling difference, we arrive instead of \eqref{bd:0} at
$$ \frac{1}{2} \d_t \big(|v(t)|_{H^s}^2\big) + |\d_x v|_{H^s}^2 \leq \left( \frac{\e t}{2\pi k_0} + C_s |v|_{L^\infty} \right) |\d_x v|_{H^s}^2,$$
and the condition on $t$ now takes the form $\e t < 2 \pi k_0.$ Back in the original time scale, this gives the uniform bound  \be \label{crude:bound} |(\e \d_x)^\a u - i t - \bar a |_{L^2} \leq |a(k_0 \cdot) -\bar a|_{H^s}.\ee
Our goal in this Section is to refine this crude a priori analysis, by proving the precise growth estimate \eqref{growth:est}.

\subsection{Approximate solution} \label{sec:app}

For the solution $v$ to \eqref{eq:vtilde}, we posit the ansatz 
\be \label{an:w} 
v = v_1 + i \e v_2, \qquad v_j \in \R,
\ee
leading to the system 
\be \label{syst:v1v2}
\big( \d_t + (\bar a + v_1) \d_x - \d_x^2 \big) \left(\begin{array}{c} v_1 \\ v_2 \end{array}\right) + \left(\begin{array}{cc} 0 & - \e^2 (t  +v_2) \\ t + v_2 & 0 \end{array}\right) \d_x  \left(\begin{array}{c} v_1 \\ v_2 \end{array}\right) = 0.
\ee
In a first step, we solve \eqref{syst:v1v2} for $\e = 0,$ corresponding to the decoupled system 
\be \label{syst:va}
 \left\{ \begin{aligned} \big( \d_t + (\bar a + v^a_1) \d_x - \d_x^2 \big) v_1^a & = 0, \\ \big( \d_t + \bar a \d_x - \d_x^2 \big) v_2^a + \d_x (v_1^a v_2^a) & = - t \d_x v_1^a,
  \end{aligned}\right.\ee
  with data 
\be \label{data:va}
 v_1^a(0,x) = a(k_0 x) - \bar a, \quad v_2^a(0) = 0.
\ee
In this Section we derive estimates for the solution to \eqref{syst:va}-\eqref{data:va}.

\subsubsection{Estimates for $v_1^a$} \label{sec:v1a}

There holds, by reality of $v_1^a,$ for $\a \leq s,$ 
$$ \big| \Re e \, (\d_x^\a (v_1^a \d_x v_1^a), \d_x^\a v_{1}^a)_{L^2} \big| \leq  |\d_x v_{1}^a|_{L^\infty} |v_1^a|_{H^s}^2 + \big| \big[ \d_x^\a, v_1^a \big] \d_x v_1^a \big|_{L^2} |v_{1}^a|_{H^s} ,$$
hence, by the Kato-Ponce inequality\footnote{We mean inequality $|[\d_x^\a, w_1] w_2|_{L^2} \leq C_s (|\d_x w_1|_{L^\infty} |w_2|_{H^{s-1}} + |w_2|_{L^\infty} |\d_x w_1|_{H^{s-1}}),$ for $0 \leq \a \leq s,$ some $C_s > 0,$ all  $w_2, \d_x w_1 \in L^\infty \cap H^{s-1}.$} 
$$ \big| \Re e \, \big( \, \d_x^\a (v_1^a \d_x v_1^a), \d_x^\a v_{1}^a \big)_{L^2} \big| \leq C_s |\d_x v_1^a|_{L^\infty} |v_1^a|_{H^s}^2,$$
for some $C_s > 0.$ 
We deduce the upper bound
\be \label{up:v1a}
 \frac{1}{2} \d_t \big( |v_1^a|_{H^{s}}^2 \big) + |\d_x v_1^a|_{H^{s}}^2 \leq C_s |\d_x v_1^a|_{L^\infty} |v_1^a|_{H^{s}}^2,
\ee
Besides, $\d_x v_1^a$ solves 
 \begin{equation} \label{03}
 (\d_t - \d_x^2 + v_1^a \d_x) (\d_x v_1^a) = - (\d_x v_1^a)^2 \leq 0,
 \end{equation}
 hence by the maximum principle
 \begin{equation} \label{04}
 |v_1^a(t)|_{L^\infty} \leq |a(k_0 \cdot) - \bar a|_{L^\infty}, \qquad |\d_x v_1^a(t)|_{L^\infty} \leq |k_0| |\d_x a|_{L^\infty}.
  \end{equation}
 Thus with \eqref{up:v1a} and Poincar\'e-Wirtinger, if $|\d_x a|_{L^\infty}$ is small enough, we obtain
 \be \label{up:v1a:2}
 \d_t \big( |v_1^a|_{H^{s}}^2 \big) + |\d_x v_1^a|_{H^{s}}^2 \leq 0,
 \ee
 implying the uniform bound
 \be \label{uniform:v1a}
 |v_1^a|_{H^{s}} \leq |a(k_0 \cdot) - \bar a|_{H^{s}}, 
 \ee 
 In particular, since $a$ is assumed to belong to $H^{s+1},$ we can repeat the argument in $H^{s+1},$ and  obtain
 \be \label{uniform:v1as+1}
 |v_1^a|_{H^{s+1}} \leq {\bm a} := |a(k_0 \cdot) - \bar a|_{H^{s+1}}.
 \ee 
  From \eqref{up:v1a:2} we deduce also
 $$ \int_0^t \t^2 |\d_x v_1^a|_{H^{s}}^2 \, d\t \leq - \int_0^t \t^2 \d_t |v_1^a(\t)|_{H^{s}}^2 \, d\t.$$
 Integrating by parts and using \eqref{uniform:v1a}, this yields
 \be \label{uniform:v1a:2}
 \int_0^t \t^2 |\d_x v_1^a(\t)|_{H^{s}}^2 \, d\t  \leq t^2 |a(k_0 \cdot) - \bar a|_{H^{s}}^2,
 \ee
 corresponding to a growth of the time moment that is better than would be expected from \eqref{uniform:v1a}.

 \subsubsection{Estimates for $v_2^a$} \label{sec:v2a} We apply $\d_x^\a$ to equation \eqref{syst:va}(ii) in $v_2^a$ and take the scalar product with $v_{2,\a}^a := \d_x^\a v_2^a,$ to obtain 
 \be \label{v2a:1} \frac{1}{2} \d_t (|v_{2,\a}^a|_{L^2}^2) + |\d_x v_{2,\a}^a|_{L^2}^2 \leq \big| t (\d_x v_{1,\a}^a, v_{2,\a})_{L^2}\big| + \big| (\d_x^{\a + 1}(v_1^a v_2^a), v_{2,\a})_{L^2}\big|.\ee
With the uniform bound \eqref{uniform:v1as+1}, there holds
   $$ \big| t (\d_x v_{1,\a}^a, v_{2,\a})_{L^2}\big| \leq t {\bm a} |v_{2,\a}^a|_{L^2},$$
where notation ${\bm a}$ is introduced in \eqref{uniform:v1as+1}, and %
$$ \big| (\d_x^{\a + 1}(v_1^a v_2^a), v_{2,\a})_{L^2} \big| \leq  {\bm a} |v_2^a|_{L^\infty} |v_{2,\a}^a|_{L^2} + \big| \big[ \d_x^{\a+1}, v_2^a \big] v_1^a \big|_{L^2} |v_{2,\a}^a|_{L^2}, \qquad \a \leq s.$$
Using Kato-Ponce again, this gives for $\a \leq s$ the bound 
$$ %
\begin{aligned} 
\frac{1}{2} \d_t (|v_{2,\a}^a|_{L^2}^2) + |\d_x v_{2,\a}^a|_{L^2}^2 & \leq  t {\bm a} |v_{2}^a|_{H^s} + {\bm a} |v_{2}^a|_{L^\infty}  |v_{2}^a|_{H^s} \\ & + C_s \big( {\bm a} |\d_x v_2^a|_{L^\infty} |v_{2}|_{H^s} + |a - \bar a|_{L^\infty} |v_2^a|_{H^{s+1}}  |v_2^a|_{H^{s}} \big).
 \end{aligned}
$$ %
If ${\bm a}$ is small enough, we can absorb the last three terms in the above upper bound by the viscous term in the left-hand side of the inequality, by use of Poincar\'e-Wirtinger. After summation over $0 \leq \a \leq s,$ we are left with 
\be \label{up:v2} \d_t (|v_{2}^a|_{H^s}^2) + |\d_x v_{2}^a|_{H^s}^2 \lesssim t |\d_x v_1^a|_{H^s} |v_2^a|_{H^s}.\ee
Integrating in time, 
 we find  
$$ |v_2^a(t)|_{H^s}^2 + \int_0^t |\d_x v_2^a(\t)|_{H^s}^2 \, d\t  \lesssim \int_0^t \t^2 |\d_x v_1^a(\t)|_{H^s}^2 \, d\t + \int_0^t |v_2^a|^2_{H^s},$$ 
and with the uniform time-integrated bound \eqref{uniform:v1a:2}, and another application of Poincar\'e-Wirtinger, we finally obtain
\be \label{v2a:final}
 |v_2^a(t)|^2_{H^s} + \int_0^t |\d_x v_2^a(\t)|_{H^s}^2 \, d\t \lesssim t^2 {\bm a}^2.
\ee
Again, the time moment has a slower growth than would be expected from \eqref{v2a:final}. Indeed, from \eqref{up:v2}, we deduce 
the bound
$$ \int_0^t \t^2 |\d_x v_2^a|_{H^s}^2 \, d\t \lesssim \int_0^t \t |v_2^a|_{H^s}^2 \, d\t +  \int_0^t \t^3 |\d_x v_1^a|_{H^s} |v_2^a|_{H^s} \, d\t.$$
With \eqref{v2a:final} and Young's inequality this gives
$$ \int_0^t \t^2 |\d_x v_2^a|_{H^s}^2 \, d\t \lesssim t^4 {\bm a}^2 + t^2 \left( \int_0^t \t^2 |\d_x v_1^a|_{H^s}^2 \, d\t + \int_0^t |v_2^a|_{H^s}^2 \, d\t\right),$$
implying with \eqref{v2a:final} and the time-integrated bound \eqref{uniform:v1a:2} the estimate 
\be \label{v2a:2}
 \int_0^t \t^2 |\d_x v_2^a|_{H^s}^2 \, d\t \lesssim t^4 {\bm a}^2.
 \ee
\subsection{The perturbation variable}

We look for the solution to \eqref{syst:v1v2} in the form
\begin{equation} \label{06}
 v_1 = v_1^a + \e w_1, \qquad v_2 = v_2^a + \e w_2,
 \end{equation}
 where $(v_1^a,v_2^a)$ is the approximate solution of Section \ref{sec:app}. Then the perturbation variable $(w_1,w_2)$ solves
\begin{equation} \label{eq:w1w2} \left\{ \begin{aligned}
 (\d_t - \d_x^2) w_1 + \d_x( v_1^a w_1) + \frac{\e}{2} \d_x \big( w_1^2 -  (t + v^a_2 + \e w_2)^2\big) & = 0, \\
 (\d_t - \d_x^2) w_2 + \d_x \big(  (v_1^a + \e w_1) w_2 \big) & = - \d_x \big( (t + v_2^a) w_1 \big).
 \end{aligned} \right.
\end{equation}
with null initial data: $w_1(0) = w_2(0) = 0.$ The equations \eqref{eq:w1w2} preserve the mean value. In particular, there holds $\bar w_1(t) = \bar w_2(t) = 0$ for all $t.$
\subsubsection{Estimates for $w_1$} \label{sec:w1} There holds
 \be \label{est:w1:0}  \frac{1}{2} \d_t (|w_{1,\a}|_{L^2}^2) + |\d_x w_{1,\a}|_{L^2}^2 + {\rm I} + \e {\rm II} + \e {\rm III} = 0, \qquad w_{1,\a} := \d_x^\a w_1.\ee
 The convective terms ${\rm I}$ and ${\rm II}$ are bounded as above, integrating by parts and using Kato-Ponce. Precisely,
 $$ | \, {\rm I} \, | := \big| \, \Re e \, (\d_x^{\a+1} (v_1^a w_1), w_{1,\a})_{L^2} \, |  \leq \big| \Re e \, (v_1^a \d_x w_{1,\a}, w_{1,\a})_{L^2}| + \big| \big[ \d_x^{\a + 1}, v_1^a \big] w_1\big|_{L^2} |w_1|_{H^s},
$$
 so that
 $$ |\, {\rm I}\, | \lesssim |\d_x v_1^a|_{L^\infty} |w_1|_{H^s}^2 + |\d_x v_1^a|_{H^{s}} |w_1|_{L^\infty} |w_1|_{H^s};
 $$
 and
 $$ |\, {\rm II}\, | := \frac{1}{2} \big|\Re e\, (\d_x^{\a+1} (w_1^2), w_{1,\a})_{L^2} \big| \lesssim |\d_x w_1|_{L^\infty} |w_1|_{H^s}^2.
 $$ 
 For the third term, we simply integrate by parts so as to let $\d_x w_{1,\a}$ appear in the $L^2$ scalar product. This gives
 $$ \frac{\e}{2} \big| \big( \d_x^\a (t + v_2)^2, w_{1,\a} \big)_{L^2} \big| \leq \e (t + C_s |v_2^a|_{H^s}) |v_2^a|_{H^s} |\d_x w_1|_{H^{s}},$$
 where we used the fact that $H^s$ is an algebra, and
 $$ \e^3 \big| \big( \d_x^{\a + 1} (w_2^2), w_{1,\a} \big)_{L^2} \big| \leq \e^3 |w_2^2|_{H^s} |\d_x w_1|_{H^{s}},$$
 and by Moser's product inequality
 $$ \e^3 \big| \big( \d_x^{\a + 1} (w_2^2), w_{1,\a} \big)_{L^2} \big| \lesssim \e^3 |w_2|_{L^\infty} |w_2|_{H^s} |\d_x w_1|_{H^{s}}.$$
Finally, for the last term in ${\rm III}$ we have
 $$ \e^2 \big| \big( \d_x^{\a+1} \big( (t + v_2^a)w_2 \big), w_{1,\a} \big)_{L^2} \big| \leq \e^2 t |w_2|_{H^s} |\d_x w_1|_{H^{s}} + \e^2 C_s |v_2^a|_{H^s} |w_2|_{H^s} |\d_x w_1|_{H^{s}}.$$ %
With the bound on $v_1^a$ from Section \ref{sec:v1a}, reliant on smallness of $|a - \bar a|_{H^{s+1}},$ the upper bound for ${\rm I}$ is absorbed by the viscous term. That is, from \eqref{est:w1:0} and the above bound on ${\rm I},$ we find
$$ \frac{1}{2} \d_t (|w_1|_{H^s}^2) + |\d_x w_1|_{H^s}^2 \leq C(a) |w_1|_{H^{s}}^2 + \e {\rm II} + \e {\rm III},$$
where $C(a) \to 0$ as ${\bm a} = |a(k_0 \cdot) - \bar a|_{H^{s+1}} \to 0,$ implying the bound
$$ \frac{1}{2} \d_t (|w_1|_{H^s}^2) + (1 - C(a)) |\d_x w_1|_{H^s}^2 \leq \e {\rm II} + \e {\rm III}.$$
We henceforth use $\s$ as a generic notation for positive constants which can be made arbitrarily small by choosing ${\bm a} = |a(k_0 \cdot) - \bar a|_{H^{s+1}}$ small. The constant $C(a)$ above falls into this category. With the above bounds on ${\rm II}$ and ${\rm III},$ we deduce the upper bound %
\be \label{dt:w1} \frac{1}{2} \d_t (|w_1|_{H^s}^2) +  (1 - \s) |\d_x w_1|_{H^s}^2  \leq A +   \e (t + C_s |v_2^a|_{H^s}) |v_2^a|_{H^s} |\d_x w_1|_{H^{s}},\ee %
 with notation
$$ %
 A := \e C_s |\d_x w_1|_{L^\infty} |w_1|_{H^s}^2 + \e^2 \big( t( 1  + C_s {\bm a}) + \e C_s |w_2|_{L^\infty} \big) |w_2|_{H^s} |\d_x w_1|_{H^{s}}.
$$ %
Using the bounds for $v_2^a$ in \eqref{v2a:final}, we find
$$ \int_0^t \e (\t + C_s |v_2^a|_{H^s}) |v_2^a|_{H^s} |w_1|_{H^{s+1}} d\t \leq \frac{C_s}{4 \eta} \e^2 t^4 {\bm a}^2 + \eta \int_0^t |\d_x w_1|_{H^{s}}^2 \, d\t, \quad \eta > 0.$$
Thus from \eqref{dt:w1} we deduce, for $\e t \leq T$ and $0 < \eta:$  
\be \label{up:w1} 
 \frac{1}{2} |w_1(t)|_{H^s}^2 + (1 - \s) \int_0^t |\d_x w_1|_{H^s}^2 \leq \frac{C_s T^2}{4 \eta} t^2 {\bm a}^2 + \int_0^t A \, d\t.
 \ee
Here $\s$ (recall our notational convention set out just above \eqref{dt:w1}) depends in particular on $\eta,$ the multiplicative constant coming from the elementary Young inequality. Going back to \eqref{dt:w1}, and using \eqref{v2a:final}, we also have 
$$ (1 - \s) \int_0^t \t^2 |\d_x w_1|_{H^s}^2 \leq 2 \int_0^t \t |w_1|_{H^s}^2 \, d\t + \int_0^t \t^2 A \, d\t + \int_0^t \e \t^3 C_s |v_2^a|_{H^s} |\d_x w_1|_{H^{s}}\, d\t.$$
The last term above is in part absorbed to the left, via
$$ \int_0^t \e \t^3 |v_2^a|_{H^s} |\d_x w_1|_{H^{s}}\, d\t \leq \frac{1}{4 \eta} \int_0^t \e^2 \t^4 |v_2^a|_{H^s}^2 \, d\t + \eta \int_0^t \t^2 |\d_x w_1|_{H^{s}}^2 \, d\t,$$
which with \eqref{v2a:2} gives
$$ \int_0^t \e \t^3 |v_2^a|_{H^s} |\d_x w_1|_{H^{s}}\, d\t \leq \frac{C_s T^2}{4 \eta} t^4 {\bm a}^2 + \eta \int_0^t \t^2 |\d_x w_1|_{H^{s}}^2 \, d\t.$$
Thus we obtained the estimate
\be \label{up:w1:2}
 (1 - \s) \int_0^t \t^2 |\d_x w_1|_{H^s}^2 \leq \frac{C_s T^2}{4 \eta}  t^4 {\bm a}^2 + 2 \int_0^t \t |w_1|_{H^s}^2 \, d\t + \int_0^t \t^2 A \, d\t.
 \ee
\subsubsection{Estimates for $w_2$} \label{sec:w2} The contribution of the linear convective term $\d_x (v_1^a w_2)$ is absorbed by the viscosity, if ${\bm a}$ (defined in \eqref{uniform:v1as+1}) is small enough, just like term ${\rm I}$ in Section \ref{sec:w1}. The nonlinear convective term $\d_x (w_1 w_2)$ is also handled as in Section \ref{sec:w1}: %
$$ \big| (\d_x^{\a+1} (w_1 w_2), w_{2,\a})_{L^2} \big| \lesssim | \d_x w_1|_{L^\infty} |w_2|_{H^s}^2 + |w_2|_{L^\infty} |w_1|_{H^{s+1}} |w_2|_{H^s}.$$
Since $H^{s+1}$ is an algebra, the source term in the right-hand side of \eqref{eq:w1w2}(ii) contributes
$$  \big| \big( \d_x^{\a+1} \big( (t + v_2^a) w_1 \big), w_{2,\a}\big)_{L^2} \big|  \leq t |\d_x w_1|_{H^{s}} |w_2|_{H^s} + C_s | v_2^a|_{H^{s+1}} |w_1|_{H^{s+1}} |w_2|_{H^s},$$
hence, with \eqref{v2a:final},
$$  \big| \big( \d_x^{\a+1} \big( (t + v_2^a) w_1 \big), w_{2,\a}\big)_{L^2} \big|  \leq t |\d_x w_1|_{H^{s}} |w_2|_{H^s} + C_s {\bm a}\big) |w_1|_{H^{s+1}} |w_2|_{H^s}.$$
Thus we obtain
\be \label{up:w2}
  \frac{1}{2} \d_t (|w_2(t)|_{H^s}^2) + (1 - C(a)) |\d_x w_2|_{H^s}^2   \leq  B,
\ee
 where $C(a) \to 0$ as ${\bm a} \to 0,$ just like in Section \ref{sec:w1}. In \eqref{up:w2}, we introduced notation
$$ B := \big( t + C_s {\bm a}\big) |\d_x w_1|_{H^s} |w_2|_{H^s} + \e C_s \big( | \d_x w_1|_{L^\infty} |w_2|_{H^s}^2 + |w_2|_{L^\infty} |w_1|_{H^{s+1}} |w_2|_{H^s}  \big).
$$ %
\subsubsection{Continuation of a priori bounds} \label{sec:cont}

The existence and uniqueness of the solution $(w_1, w_2)$ to \eqref{eq:w1w2} issued from $(0,0)$ is granted in very small time\footnote{For instance by the arguments of Section \ref{sec:crude}.}, with bounds 
\begin{equation}\label{hyp:w12}
|\d_x w_1|_{L^\infty} \leq M_1 (1 + t), \quad |w_2|_{L^\infty} \leq M_2 (1 + t^2),
\end{equation}
We now show that the a priori estimates of Sections \ref{sec:w1} and \ref{sec:w2} imply that bounds \eqref{hyp:w12} are propagated up to time $T/\e,$ where $T$ can be arbitrarily close to $2 \pi k_0,$ if $|a - \bar a|_{H^{s+1}}$ is small enough.  

So long as \eqref{hyp:w12} holds, and $\e t \leq T,$ we may bound $A$ and $B$ by\footnote{For $t \geq 1,$ corresponding to $t \geq \e$ in the original time scale, up to changing $M_1$ and $M_2$ into $M_1/2$ and $M_2/2.$} 
\be \label{AB} \begin{aligned} A & \leq  C_s M_1\, T |w_1|^2_{H^s} + \e^2 t (1 + \s)  |w_2|_{H^s}|\d_x w_1|_{H^{s}}, \\
B & \leq C_s M_1 \, T |w_2|^2_{H^s} + t (1 + \s) |w_2|_{H^s}|\d_x w_1|_{H^{s}}.\end{aligned}\ee
where $\s = C_s |a - \bar a|_{H^{s+1}} + C_s M_2 T,$ consistent with the notational convention that we set out just above \eqref{dt:w1} if we allow $\s$ to depend on $M_2.$ 
Plugging in \eqref{up:w1} and \eqref{up:w2}, this gives, for $w = (w_1,w_2),$ the bound 
\be \label{for:up0} \frac{1}{2} |w|^2_{H^s} + (1 - \s) \int_0^t |\d_x w|^2_{H^s} \leq \frac{C_s T^2}{4 \eta} t^2 |a - \bar a|_{H^{s+1}}^2 + \int_0^t  C_s M_1 T  |w|^2_{H^s} d \t + D,
\ee
with notation
\be \label{D} D := (1 + \s) \int_0^t \t |w_2|_{H^s} |\d_x w_1|_{H^{s}} \, d\t,\ee
where $\s$ depends in particular on $M_1.$ %
If $M_1$ is small enough, depending on $T$ and $s,$ we may absorb, via Poincar\'e-Wirtinger, the second term in the upper bound of \eqref{for:up0} in the viscous term in the left-hand side, and obtain 
\be \label{for:up1}
 \frac{1}{2} |w|^2_{H^s} + (1 - \s) \int_0^t |\d_x w|^2_{H^s} \leq \frac{C_s T^2}{4 \eta} t^2 {\bm a}^2 + D
\ee
We now concentrate on a bound for $D.$ We are going to absorb most terms in $D$ into the viscous term in the left-hand side. To this end we will use the time-integrated bound \eqref{up:w1:2} for $w_1.$  By Young's inequality,  
\be \label{for:up} \begin{aligned} (1 + \s)^{-1} D & \leq \frac{1}{4 \eta_1}\int_0^t |w_2|_{H^s}^2 \, d\t + 
\eta_1 \int_0^t \t^2 |\d_x w_1|_{H^{s}}^2 \, d\t, \quad 0 < \eta_1.\end{aligned}
\ee
 Under \eqref{hyp:w12}, by \eqref{AB}(i), for $0 < \eta_2:$
 $$ \begin{aligned} \int_0^t \t^2 A \, d\t & \leq C_s M_1 \int_0^t \t^2 |w_1|_{H^s}^2 \, d\t + (1 + \s) \frac{(\e t)^2}{4 \eta_2} \int_0^t |w_2|_{H^s}^2 \, d\t \\ & \quad + (1 + \s) \eta_2 (\e t)^2 \int_0^t \t^2 |\d_x w_1|_{H^s}^2 \, d\t.\end{aligned}$$
Using the above bound in \eqref{up:w1:2}, we obtain 
\be \label{E} %
\Big(1 - \s - (1 + \s) \eta_2 (\e t)^2\,\Big) \int_0^t \t^2 |\d_x w_1|_{H^s}^2 \, d\t \leq (1 + \s) \frac{(\e t)^2}{4 \eta_2} \int_0^t |w_2|_{H^s}^2 + E,
\ee
where 
$$ E := \frac{C_s}{4 \eta}  t^4 {\bm a}^2 + 2 \int_0^t \t |w_1|_{H^s}^2 \, d\t + C_s M_1 \int_0^t \t^2 |w_1|_{H^s}^2.$$
We may overlook the last term in $E,$ by Poincar\'e-Wirtinger, up to changing $\s$ in \eqref{E}. For the second term in $E,$ we use \eqref{up:w1}:
$$ \int_0^t \t |w_1|_{H^s}^2 \, d\t \lesssim t \int_0^t |\d_x w_1|_{H^s}^2 \lesssim \frac{C_s}{\eta} t^3 {\bm a}^2 + t \int_0^t A \, d\t.$$
We use \eqref{AB}(i) again: the contribution of the first term in $A$ is absorbed by the left-hand side, if $M_1$ is small enough. The contribution of the second term in $A$ has the same form as $D,$ with an extra $\e^2 t$ multiplicative prefactor. Hence we can go through the same steps, and obtain \eqref{E} where $E$ consists only of its first term, up to changing the values of $\s.$ %

Back to \eqref{for:up}, we obtained
$$  D \leq K \int_0^t |w_2|_{H^s}^2 + \frac{C_s}{4 \eta}  t^4 {\bm a}^2,$$
where
$$ K := (1 + \s) \Big( \frac{1}{4 \eta_1} + \eta_1 (1 + \s) \big(1 - \s - \eta_2 (\e t)^2\big)^{-1} \frac{(\e t)^2}{4 \eta_2}\Big),$$
which we may write, with our notational convention for $\s,$ 
$$ K = \s + \frac{1}{4 \eta_1} + \eta_1 \big(1 - \eta_2 (\e t)^2\big)^{-1} \frac{(\e t)^2}{4 \eta_2}.$$
There holds
$$ \min_{\eta_1 > 0} K = \s + \big(1 - \eta_2 (\e t)^2\big)^{-1/2} \frac{(\e t)}{2 \eta_2^{1/2}}.$$
By Poincar\'e-Wirtinger \eqref{pw2}, since $w_2$ has zero mean and depends on $x$ through $k_0 x,$ 
$$ \int_0^t |w_2|_{H^s}^2 \, d\t \leq \frac{1}{(2 \pi k_0)^2} \int_0^t |\d_x w_2|_{H^s}^2 \, d\t.$$
Thus we can absorb the first term in $D$ by the left-hand side of \eqref{for:up1} for $t$ and $\eta_2$ such that
\be \label{cond:t:eta2}
 \big(1 - \eta_2 (\e t)^2\big)^{-1/2} \frac{\e t }{2 \eta_2^{1/2}} < (2 \pi k_0)^2,
\ee
and since 
$$ \min_{\eta_2 > 0} \big(1 - \eta_2 (\e t)^2\big)^{-1/2} \frac{\e t }{2 \eta_2^{1/2}} = (\e t)^2,$$
we see that \eqref{cond:t:eta2} amounts to $\e t < 2 \pi k_0,$ corresponding to a limiting time arbitrarily close to the transition time.

Thus we proved that in \eqref{for:up1} all the terms in $D$ save for the first term in $E$ can be absorbed in the left-hand side. This gives 
\be \label{for:up2}
  \frac{1}{2} |w|^2_{H^s} + (1 - \s) \int_0^t |\d_x w|^2_{H^s} \leq \frac{C_s}{\eta} t^4 {\bm a}^2,
 \ee
 under the conditions that in \eqref{hyp:w12} the constants $M_i$ be small enough, and for $\eta = \eta({\bm a}, T, M_i),$ where $\eta \to 0$ as $T \to 2 \pi k_0.$ 
We will use \eqref{for:up2} as a bound in $w_2,$ and now go back to the bound \eqref{up:w1} in $w_1.$ By \eqref{AB}(i), for $\s \leq 1$ there holds 
$$ %
 \int_0^t A \leq  C_s M_1 T \int_0^t  |w_1|^2_{H^s} \, d\t  +  \frac{\e^2}{\eta'} \int_0^t |w_2|^2_{H^s} + \eta' (\e t)^2\int_0^t |\d_x w_1|^2_{H^{s}}\,,
$$ %
for any $\eta' > 0,$ and now using \eqref{for:up2},
$$ %
 \int_0^t A \, d\t \leq \Big( C_s M_1 T + \frac{\eta' (\e t)^2}{(2 \pi k_0)^2}\Big) \int_0^t |\d_x w_1|_{H^s}^2 \, d\t + \frac{C_s (\e T)^2}{\eta \eta'} t^4 {\bm a}^2\,.
$$ %
With the above bound and \eqref{up:w1}, we finally obtain, if $M_1$ and $\eta'$ are small enough, 
\be \label{final:w1}
 \frac{1}{2} |w_1|^2_{H^s} + (1 - \s) \int_0^t |\d_x w_1|^2_{H^s} \leq \frac{C_s T^2}{\eta}\Big(1 + \frac{1}{\eta'} \Big) t^2 {\bm a}^2.
\ee

From \eqref{for:up2} and \eqref{final:w1}, we conclude that the a priori bounds \eqref{hyp:w12} propagate up to $T/\e,$ as follows. 

Given $0 < \e$ and $T < 2 \pi k_0,$ the constraint $\s < 1$ takes the form
 \be \label{cons:sigma}
 C(a) + C_s(M_1 + M_2) + \eta + \eta' < (1 - \e) \Big(1 - \frac{T}{2 \pi k_0}\Big),
 \ee
 where $0 < C(a) \to 0$ as ${\bm a} \to 0.$ We first choose $M_i$ and $\eta, \eta',$ such that $C_s( M_1 + M_2) + \eta + \eta'$ equals one half of the upper bound in \eqref{cons:sigma}. This means in particular $\eta,\eta' \to 0$ as $T \to 2 \pi k_0.$ Then we choose ${\bm a}$ small enough so that not only \eqref{cons:sigma} is satisfied, but also
 \be \label{cons:a}
 {\bm a} \lesssim \eta (M_1^2 +  M_2^2),
 \ee
 the implicit multiplicative constant depending only on $s$ and $k_0.$ Then \eqref{for:up2} and \eqref{final:w1} imply that the pointwise a priori bounds \eqref{hyp:w12} propagate up to $T/\e.$ 
 
 We obtained closed a priori bounds in $C^0([0, T/\e], H^s).$ This translates into an existence and uniqueness result by classical arguments. For the original variable $u:$ 
$$ u = \bar a + v_1^a + \e w_1 + i \big( t + \e v_2^a + \e^2 w_2\big),$$
the bounds \eqref{uniform:v1a} and \eqref{v2a:final} on the approximate solution $(v_1^a, v_2^a)$ and \eqref{for:up2} and \eqref{final:w1} on the perturbation $(w_1,w_2)$ translate into estimates \eqref{growth:est}, valid for ${\bm a} = |a(k_0 \cdot) - \bar a|_{H^{s+1}}$ small enough, depending on $1 - T/(2 \pi k_0).$

\section{Proof of Theorem \ref{th2}: small oscillating data and growth in time $O(1)$} \label{sec:growth}

 We consider here small, highly oscillatory data $\dsp{u(0,x)  = \e^{\a} a(k_0 x/\e)}$
 with $\a > 1/3$ and $a \in H^1(\T),$ such that $\dsp{a_1 := \int_{\T} e^{-2 i \pi x} a(x) \, dx \neq 0.}$

\subsection{Uniform bounds in short time} \label{sec:small}

In a first step, we prove a short-time existence result. We could simply use the result of Section \ref{sec:crude} here, but the notation and estimates of the present Section will be useful later on. We posit the ansatz
\be \label{an}
 u(\e,t,x) = i t + \bar a + \e^\a v\left(\frac{t}{\e}, \frac{x}{\e}\right).
\ee
Then $v$ solves
\be \label{eq:vtildebis} %
\d_t v + (i \e t + \bar a + \e^\a v) \d_x v -  \d_x^2 v  = 0, \qquad \bar v \equiv 0, \qquad v(\e,0,x) = a(k_0 x) - \bar a.
\ee
As in Section \ref{sec:proof1}, so long as it is defined $v(t)$ is a function of $k_0 x.$ We denote $(v_k)$ the Fourier mode corresponding to the $k$-th harmonics of $k_0:$ 
\be \label{def:v1} v(\e,t,x) = \sum_{k \in \Z} e^{2 i \pi k \cdot k_0 x} v_k(\e,t).\ee
There holds $v_1(0,x) = a_1 \neq 0.$ The leading harmonics will be $v_1,$ in the sense that $v_1$ will grow from $t = 2 \pi k_0$ on, before all other harmonics. The goal is to reach the amplification time $4 \pi k_0,$ in the original time scale. In the time scale for $v,$ this corresponds to $4 \pi k_0/\e.$ We let
$$w(t,x) := v(t,x) - e^{ i x} v_1(t) = \sum_{k \notin \{0,1\}} e^{2 i \pi k \cdot k_0 x} v_k(t) =: \sum_{k \in \Z} e^{2 i \pi  k \cdot k_0x} w_k(t).$$
Then, $v_1$ and $w$ solve the triangular system
\be \label{syst:v1w}
 \left\{\begin{aligned}
  \d_t v_1 + 2 \pi k_0\big( \, i \bar a + 2 \pi k_0 - \e t \big) v_1 & = \e^\a f(w,w), \\
  \d_t w + (i \e t + \bar a) \d_x w - \d_x^2 w & = \e^\a \d_x g(v,v),
  \end{aligned}\right.
  \ee
 with notation
\be \label{def:f-g}
 f(w,w) :=  - 2 i \pi k_0 \sum_{k_1 + k_2 = 1} w_{k_1} w_{k_2}, \quad g(v,v) := - \sum_{k \neq 1} e^{2 i \pi k_0 \cdot k x} \sum_{k_1 + k_2 = k} v_{k_1} v_{k_2}.
\ee
In integral form, 
\be \label{syst:v1wbis} \begin{aligned}
 v_1(t) & = e^{\l_1(t)} v_1(0) + \e^\a \int_0^t e^{\l_1(t) - \l_1(\t)} f(w,w)(\t) \, d\t, \\
 w_k(t) & = e^{\l_k(t)} w_k(0) + \e^\a (2 i \pi k_0) \, k  \int_0^t e^{\l_k(t) - \l_k(\t)} \widehat{g(v,v)}(\t,k) \, d\t.
 \end{aligned}\ee
where $\l_k(t)$ stands for
\be \label{def:lambda}
\l_k(t) := - 2 i \pi k_0 k \bar a t + \pi k_0 k \big(\e t^2 -  4 \pi k_0 k t \big).
\ee
   There holds, by Young's convolution inequality and Parseval's equality, the pointwise bounds
\be \label{young}
 |f(z_1,z_2)| \leq |z_1|_{L^2} |z_2|_{L^2}, \quad |\widehat{g(z_1,z_2)}| \leq |z_1|_{L^2} |z_2|_{L^2}.
\ee
For $0 \leq t$ and $k \leq -1,$ we will use the elementary bounds 
 $$ e^{\re \l_k(t)} \leq e^{- t (2 \pi k_0 k)^2}, \quad \int_0^t e^{\re \l_k(t) - \re \l_k(\t)} |k| \, d\t \leq \int_0^t e^{- (t - \t) (2 \pi k_0 k)^2} |k| \, d\t \lesssim \frac{1}{|k|}.$$
For $3 \leq k$ and $\dsp{0 \leq \t \leq t \leq \frac{5 \pi k_0}{\e}},$  %
 $$ \re \l_k(t) = - t (2 \pi k_0)^2 k (k-5/4),$$
 and
 \be \label{elem} \begin{aligned} \re \l_k(t) - \re \l_k(\t) &  = (t - \t) \pi k_0 k  \Big( \e(t + \t) - 4 \pi k_0 k \Big) \\ & \leq - (t  - \t) (2 \pi k_0)^2 k (k-5/2).\end{aligned}\ee
We deduce the a priori bound, for $0 \leq t \leq 5 \pi k_0/\e,$ where $5 \pi k_0/\e$ is only a convenient limiting time that is greater than our target $4 \pi k_0/\e:$ 
\be \label{ap:w}
 |w_k(t)| \lesssim e^{- t (2 \pi k_0)^2 |k|} |w_k(0)| + \frac{\e^\a}{|k|} |v|_{L^\infty([0,t],L^2(\T))}^2\,, \quad k \leq -1 \,\, \,\mbox{and}\,\,\, 3 \leq k.
\ee
From the equality in \eqref{elem}, we see that we cannot obtain a good decay estimate for the Duhamel term in $w_2$ up to $t = 4 \pi k_0/\e.$ For $0 \leq \t \leq t \leq 2 \pi k_0/\e,$ 
\be \label{ap:w20}
 |w_2(t)| \lesssim e^{- t(2 \pi k_0)^2 k } |w_2(0)| + \e^\a |v|_{L^\infty([0,t],L^2(\T))}^2.
\ee
By summation in $\ell_2,$ the bounds \eqref{ap:w} and \eqref{ap:w20} give %
\be \label{ap:w2}
 |w(t)|_{L^2} \lesssim e^{-t (2 \pi k_0)^2} |a(k_0 \cdot) - \bar a|_{L^2} + \e^\a |v|_{L^\infty([0,t], L^2(\T))}^2\,. 
 \ee
In order to bound the equation in $v_1,$ we compute
$$ \int_0^t e^{\re \l_1(t) - \re \l_1(\t)} \, d\t  = e^{\pi k_0 \e (t - 2 \pi k_0/\e)^2} \int_{0}^{t} e^{- \pi k_0 \e (\t - 2 \pi k_0/\e)^2} \, d\t.$$
This last error function integral can be evaluated for $t < 2 \pi k_0/\e:$ indeed we have then 
\be \label{errfn} \begin{aligned} \int_{0}^{t} e^{- \pi k_0 \e (\t - 2 \pi k_0/\e)^2} \, d\t  & = (\pi k_0 \e)^{-1/2} \int_{(\pi k_0 \e)^{1/2} (2 \pi k_0/\e -t)}^{2 (\pi k_0)^{3/2}\e^{-1/2}} e^{-z^2}\, dz  \\ & \leq (\pi k_0 \e)^{-1/2}  \int_{(\pi k_ 0\e)^{1/2}(2 \pi k_0/\e - t)}^{\infty} e^{-z^2} \, d\z \\ & \leq \frac{1}{2 \pi k_0 (2 \pi k_0 - \e t)} e^{- \pi k_0 \e (t - 2 \pi k_0/\e)^2}, \end{aligned}
\ee
where in the last inequality we used
$$ \int_x^\infty e^{-z^2} \, dz \leq \frac{e^{-x^2}}{2x}, \qquad 0 < x.$$
This gives, for $0 \leq t < 2 \pi k_0/\e:$ 
\be \label{ap:v} %
 |v_1(t)| \lesssim |v_1(0)| + \frac{\e^\a}{2 \pi k_0 - \e t} |w|_{L^\infty([0,t],L^2(\T))}^2.
\ee %
The a priori bounds, \eqref{ap:w2} and \eqref{ap:v}, are uniform with respect to $\e$ over times intervals $[0, T(\e)/\e],$ with $T(\e) = 2 \pi k_0 - O(\e^\a).$ 
\subsection{Beyond the transition time $t =2 \pi k_0$}  \label{sec:small:growth}

 In the estimates of the previous Section, the main error came from the error function estimate \eqref{errfn} associated with the propagator in the equation in $v_1.$ In the current Section, we factor out this propagator. This allows for sharper estimates. We let, using notation $\l_1$ introduced in \eqref{def:lambda}, 
 $$ \tilde v_1(t) := e^{- \re \l_1(t)} v_1(t), \qquad \tilde w_k(t) := e^{- \re \l_1(t)} w_k(t), \quad \tilde w(t,x) := \sum_{k \notin \{0,1\}} e^{2 i \pi k_0 \cdot k x} \tilde w_k(t).$$
We also denote $\tilde v(t) := e^{- \re \l_1(t)} v(t),$ with Fourier coefficients $\tilde v_k$ (equal to $\tilde w_k$ if $k \neq 1,$ to $\tilde v_1$ otherwise). Then, $\tilde v_1$ and $\tilde w_k$ solve 
$$ %
 \begin{aligned}
 \tilde v_1(t) & = v_1(0) + \e^\a \int_0^t f(\tilde w, w)(\t) \, d\t,  \\
 \tilde w_k(t) & = e^{\mu_k(t)} w_k(0) + \e^\a (2 i \pi k_0) k  \int_0^t e^{\mu_k(t) - \mu_k(\t)} \widehat{g(\tilde v, v)}(\t,k)\,d\t,
  \end{aligned} $$ %
  with notation
 \be \label{def:mu}
 \mu_k(t) := \l_k(t) - \re \l_1(t) = - 2 i \pi k_0 k \bar a t + \pi k_0 (k-1) \e t^2 - 4 (\pi k_0)^2 (k^2 - 1) t,
 \ee 
 so that
 \be \label{re:mu}
 \re \mu_k(t) = \pi k_0 (k-1) t \big( \e t - 4 \pi k_0 (k+1) \big).
 \ee
Let $M \geq 2$ and $t_\star(\e,M)$ such that  
\be \label{unif} |v(t) | \leq M | a(k_0 \cdot) - \bar a |_{L^2}.
\ee
uniformly in $t \in [0, t_\star(\e)]$ and in $\e \in (0,\e_0),$ for $\e_0$ small enough (depending on $M$), where $t_\star(\e,M)$ is a final observation time, smaller than $5 \pi k_0/\e.$ We denote $t_\star(\e) = t_\star(\e,M)$ in the following.
The analysis of the previous Section implies that such a bound exists if $t_\star(\e) = T(\e)/\e = 2 \pi k_0/\e - O(\e^{\a-1}).$ Our goal here is to extend this limiting time.

\subsubsection{Large Fourier indices} \label{sec:neg} For $k \leq -2,$ the corresponding modes are decaying, since $\re \mu_k(t)$ is negative and decaying for all $t > 0.$  Here we can simply overlook the contribution of $i t \d_x$ to $\re \mu_k,$ and obtain, for $t \leq t_\star(\e):$ 
\be \label{est:tildew-2}
|\tilde w_k(t)| \lesssim e^{- 4 \pi k_0 t (k^2 - 1)} |w_k(0)| + \frac{\e^\a M}{|k|} |\tilde v|_{L^\infty L^2}, \qquad k \leq -2. 
\ee 
The modes associated with $k \geq 2$ are not growing in our observation window $[0, 5 \pi k_0/\e].$ Indeed, for $0 \leq \t \leq t \leq 5 \pi k_0/\e$ and $2 \leq k,$ there holds, by \eqref{re:mu}, 
$$ \re \mu_k(t) = t (k-1) \Big( \frac{\e t}{2} - (k+1) \Big) \leq - (2 \pi k_0)^2 t (k^2 -1/4),$$
and
$$ \begin{aligned} \re \mu_k(t) - \re \mu_k(\t) & = \pi k_0 (t - \t) (k-1) \big( \, \e (t + \t)  - 4 \pi k_0 (k+1) \,\big) \\ & \leq -  (2 \pi k_0)^2 (t - \t) (k^2 - 3/2).\end{aligned}$$
Thus for $t \leq t_\star(\e):$ 
\be \label{est:tildewk2}
|\tilde w_k(t)| 
\lesssim e^{- (2 \pi k_0)^2 t (k^2 - 1/4)} |w_k(0)| + \frac{\e^\a M}{|k|} |\tilde v|_{L^\infty L^2}, \quad 2 \leq k.
\ee

\subsubsection{Small Fourier indices} These are $k = -1$ and $k = 1.$ For $k = -1,$ the diffusion fails to provide decay. We could simply use the convolution bound \eqref{young} and obtain
\be \label{w-1} |\tilde w_{-1}(t)| \lesssim e^{- 2 \pi k_0 \e t^2} |w_{-1}(0)| + \e^\a M |\tilde v|_{L^\infty L^2} \int_0^t e^{-2 \pi k_0 \e (t^2 - \t^2)} \, d\t,\ee 
but this would bring out a factor $\e^{\a - 1/2},$ as in \eqref{errf2} below, and impose condition $\a > 1/2.$ The bound \eqref{w-1} can be refined as follows. 
 
We observe that, since the mean mode is identically zero, all the terms in the sum $\widehat{g(\tilde v, v)}(t,-1)$ have the form $\tilde w_{k_1} v_{k_2},$ with $k_1 \leq -2.$ For $v_{k_2}$ we use the postulated bound \eqref{unif}. For $\tilde w_{k_1}$ we use the decaying bound \eqref{est:tildew-2}, so that the contribution of $\tilde w_{k_1} v_{k_2}$ to the Duhamel term in $\tilde w_{-1}$ is controlled by 
\be \label{1st:-1} \e^\a M \int_0^t e^{-2 \pi k_0 \e (t^2 - \t^2)} \Big(e^{- 4 \pi k_0 \t (k_1^2 - 1)} |w_{k_1}(0)| + \frac{\e^\a M}{|k_1|} |\tilde v|_{L^\infty L^2} \Big) \, d\t,
\ee
 up to a multiplicative constant, for $t \leq t_\star(\e),$ with $k_1 \leq -2.$ The first term above is controlled by $\e^\a M k_1^{-2} |w_{k_1}(0)|_{L^2}.$ For the second term in \eqref{1st:-1}, we use the elementary bound
 \be \label{errf2}
 \begin{aligned} \int_0^t e^{- 2 \pi k_0 \e (t^2 - \t^2)} \, d\t & = e^{-2 \pi k_0 \e t^2} (2 \pi k_0\e)^{-1/2} \int_0^{(2 \pi k_0 \e)^{1/2} t} e^{z^2} dz \\ & \leq e^{-2 \pi k_0 \e t^2} (2 \pi k_0\e)^{-1/2}  \Big( \int_0^1 e^{z^2} \, dz + \int_1^{(2 \pi k_0 \e)^{1/2} t} z e^{z^2} \, dz\Big) \lesssim \e^{-1/2}.
\end{aligned}
\ee 
Summing over $k_1,$ we obtain 
\be \label{est:tildew-1}
 |\tilde w_{-1}(t)| \lesssim e^{-(2 \pi k_0) \e t^2} |w_{-1}(0)| + \e^\a M |v(0)|_{L^2} + \e^{2 \a - 1/2}M^2 |\tilde v|_{L^\infty L^2}.
 \ee

We turn to $k = 1.$ All the terms in $f(\tilde w, w)$ have the form $\tilde w_{k_1} w_{k_2},$ where $k_1 \geq 2$ and $k_2 \leq -1.$ For $\tilde w_{k_1}$ we use \eqref{est:tildewk2}. For $w_{k_2}$ we use \eqref{ap:w}, and obtain a control by %
\be \label{for:v1} \begin{aligned}  \e^\a \int_0^t  \Big( e^{- (2 \pi k_0)^2 \t (k_1^2 - 1/4)} |w_{k_1}(0)| & + \frac{\e^\a M}{|k_1|} |\tilde v|_{L^\infty L^2} \Big) \\ & \times \Big(e^{- (2 \pi k_0)^2 \t |k_2|} |w_{k_2}(0)| + \frac{\e^\a M^2}{|k_2|} \Big) \, d\t.\end{aligned}\ee
Bounding from above the time integrals and summing over $k_1,$ we find that the above term is controlled by 
\be \label{for:ap:tildev1}
 \e^\a |v(0)|_{L^2}^2 (1 + \e^\a |\tilde v|_{L^\infty L^2}) +  t \e^{3 \a} M^3 |\tilde v|_{L^\infty L^2}.
\ee
This implies for $\tilde v_1$ the bound, for $t \leq t_\star(\e),$ 
\be \label{ap:tildev1}
|\tilde v_1(t)| \lesssim |v_1(0)| + 
 \e^\a |v(0)|_{L^2}^2 (1 + \e^\a |\tilde v|_{L^\infty L^2}) +  t \e^{3 \a} M^3 |\tilde v|_{L^\infty L^2}.
\ee

\subsubsection{Continuation of a priori bounds} \label{sec:cont:growth} We now gather \eqref{est:tildew-2}, \eqref{est:tildewk2}, \eqref{est:tildew-1} and \eqref{ap:tildev1}. Summing over $k$ and taking into account the definition of $\tilde v$ at the beginning of Section \ref{sec:small:growth}, we obtain for $\e$ small enough (depending on $M$) the a priori bound
\be \label{ap:tildev}
 | v(t)|_{L^2} \leq 2 e^{\re \l_1(t)} |v(0)|_{L^2}, \qquad t \leq \min t_\star(\e),
\ee
By definition,
$$ \re \l_1(t) = \e \pi k_0 \Big( t - \frac{2 \pi k_0}{\e} \Big)^2 - \frac{4 (\pi k_0)^3}{\e}.$$
We now let $\dsp{t_\star(\e) = \frac{4 \pi k_0}{\e} + T_\star.}$ So long as 
$$ 
 (2 \pi k_0)^2 T_\star + \e \pi k_0 T^2 \leq \ln\left(\frac{M}{2}\right),
$$ 
we see from \eqref{ap:tildev} that the a priori bound \eqref{unif} propagates beyond the transition time $t = 2 \pi k_0/\e,$ beyond the amplification time $t = 4 \pi k_0/\e,$ and up to $t_\star(\e)$ as defined above.  

\subsection{Amplification} \label{sec:exp} 

Integrating the equation in $\tilde v_1$ over $[0,t],$ with $t \leq t_\star(\e),$ the observation time $t_\star(\e)$ being defined just above, we find by the triangular inequality the lower bound  
$$ |\tilde v_1(t)| \geq |v_1(0)| - \e^\a \int_{0}^{t_\star(\e)} |f(\tilde w, w)(\t)| \, d\t.$$
We can bound the above Duhamel term with \eqref{for:ap:tildev1}, and this gives
$$ |\tilde v_1(t)| \geq |v_1(0)| -  \e^\a |v(0)|_{L^2}^2 (1 + \e^\a |\tilde v|_{L^\infty L^2}) + \e^{3 \a-1} M^3 |\tilde v|_{L^\infty L^2}.$$
In particular, for $\e$ small enough, 
$$ |\tilde v_1(t)| \geq \frac{1}{2} |v_1(0)|.$$
Back in the original variable, this gives the exponential lower bound 
\be \label{low:exp} |v_1(t)| \geq \frac{1}{2} e^{\re \l_1(t)} |v_1(0)| = \frac{1}{2} e^{\e \pi k_0 ( t - 2 \pi k_0\e )^2 - 4 (\pi k_0)^3/\e} |v_1(0)|,
\ee
as claimed in the statement of Theorem \ref{th2}. Indeed, by definition of $v$ in \eqref{an} and $v_1$ in \eqref{def:v1}, there holds $\dsp{v_1(t) = \e^{-1} \int_{\T} e^{-2 i \pi k_0 x/\e} u(\e t, x) \, dx.}$ 
In particular, evaluating at $t = t_\star(\e),$ where the final observation time is defined in Section \ref{sec:cont:growth} above, we see that if $M$ is large enough (implying $\e$ small enough), then $|v_1(t_\star(\e))| \geq 4 |v(0)|_{L^\infty},$ implying if $\bar a = 0$ the lower bound $|u(\e t_\star(\e))|_{L^\infty} \geq 4 |u(0)|_{L^\infty},$ as in the final time-step of our simulations. 

{\footnotesize \begin{thebibliography}{longueurmax}

%
%


%

\bibitem{KK} G. Kreiss, H. Kreiss, {\it Convergence to steady state of solutions of
Burgers' equation}, Appl.
Numer. Math. 2 (1986), no. 3-5, 161-179.
 \bibitem{Le} N. Lerner, {\it Metrics on the Phase Space and Non-Selfadjoint Pseudodifferential Operators,} Pseudo-Differential Operators. Theory and Applications, 3. Birkh\"auser 2010. xii+397 pp. 
 \bibitem{LMX} N. Lerner, Y. Morimoto, C. J. Xu, {\it Instability of the Cauchy-Kovalevskaya solution for a class of non-linear systems}, American J. Math., 132 (2010), no 1, 99-123.
\bibitem{LNT} N. Lerner, T. Nguyen, B. Texier, {\it The onset of instability in first-order systems}, {\tt arXiv:1504.04477}.
\bibitem{em4} L.~Yong, B.~Texier, {\it A stability criterion for high-frequency oscillations}, {\tt arXiv:1307.4196}, to appear in M\'em. Soc. Math. Fr. 
\bibitem{Me} G. M\'etivier, {\it Remarks on the well-posedness of the nonlinear Cauchy problem}, Geometric analysis of 
PDE and several complex variables, Contemp. Math., vol. 368, Amer. Math. Soc., Providence, 
RI, 2005, 337-356.
%
%
\bibitem{meta} C. Mascia, M. Strani, {\it Metastability for nonlinear parabolic equations with
applications to viscous scalar conservation laws,} SIAM J. Math. Anal. 45
(2013), no 5, 3084--3113.
\bibitem{Mo} B. Morisse, {\it On hyperbolicity and well-posedness in Gevrey spaces}, in preparation.
%
\bibitem{T4} B. Texier, {\it Approximations of pseudo-differential flows}, {\tt arXiv:1402.6868}, to appear in Indiana Univ. Math. J.
\bibitem{Z} M. Zworski, {\it Semiclassical analysis}, Graduate Studies in Mathematics, 138. American Mathematical Society, 2012. xii+431 pp.


\end{thebibliography}}
\end{document}